\pdfoutput=1
\RequirePackage{scrlfile}
\BeforePackage{cleveref}{\RequirePackage[algo2e,linesnumbered,ruled,vlined]{algorithm2e}}
\PassOptionsToPackage{compress}{cleveref}
\PreventPackageFromLoading{hypcap}
\documentclass[final,onefignum,onetabnum,onealgnum]{siamart220329}
\UnPreventPackageFromLoading{hypcap}
\usepackage[figure,figure*,table*]{hypcap}

\usepackage[T1]{fontenc}
\usepackage[utf8]{inputenc}

\usepackage{todonotes}
\usepackage{afterpage}

\usepackage{amssymb}
\usepackage{amsfonts}
\usepackage{mathtools}
\usepackage{bm}

\makeatletter
\DeclareSymbolFont{extraup}{U}{zavm}{m}{n}
\DeclareMathSymbol{\newcheckm@rk}{\mathalpha}{extraup}{128}
\DeclareMathSymbol{\newcrossm@rk}{\mathalpha}{extraup}{129}
\newcommand{\newcheckmark}{\scalebox{1.2}{$\newcheckm@rk$}}
\newcommand{\newcrossmark}{\scalebox{0.9}{$\newcrossm@rk$}}
\makeatother

\usepackage{multirow}
\usepackage{array}
\usepackage{adjustbox}

\usepackage{afterpage}
\usepackage{ragged2e}

\usepackage{etoolbox}

\makeatletter
  \def\@algocf@pre@ruled{}
  \patchcmd{\algocf@captiontext}{\endgraf}{}{}{}
\makeatother

\SetAlgoCaptionLayout{MyAlgoCaptionLayout}
\SetAlCapSty{MyAlCapSty}
\SetAlCapNameSty{MyAlCapNameSty}
\SetAlCapFnt{\small\scshape}
\SetAlCapNameFnt{\small\itshape}
\SetAlFnt{%
  \small%
  \setlength\abovedisplayskip{5pt plus 2pt minus 2.5pt}%
  \setlength\belowdisplayskip{5pt plus 2pt minus 2.5pt}%
  \setlength\abovedisplayshortskip{0pt plus 1pt}%
  \setlength\belowdisplayshortskip{2pt plus 1pt minus 1pt}%
}
\SetAlgoCaptionSeparator{}
\SetAlgoHangIndent{0pt}
\DontPrintSemicolon
\let\oldnl\nl
\newcommand{\nonl}{\renewcommand{\nl}{\let\nl\oldnl}}

\creflabelformat{line}{#2\NlSty{#1}#3}
\crefname{line}{line}{lines}
\Crefname{line}{Line}{Lines}
\crefmultiformat{subequation}%
  {\edef\crefstripprefixinfo{#1}\textup{(#2#1#3}}%
  {, \textup{#2\crefstripprefix{\crefstripprefixinfo}{#1}#3)}}%
  {, \textup{#2\crefstripprefix{\crefstripprefixinfo}{#1}#3}}%
  {, and~\textup{#2\crefstripprefix{\crefstripprefixinfo}{#1}#3)}}

\usepackage{fnpct}

\usepackage{pgfplots}
\pgfplotsset{compat=1.17}
\newlength\figheight
\newlength\figwidth

\DeclarePairedDelimiter\norm{\|}{\|}
\DeclarePairedDelimiter\snorm{\lvert}{\rvert}
\makeatletter
\newcommand\doraisetag[1]{\global\shifttag@true\raisetag{#1}}
\makeatother

\makeatletter
\if@review
  \newcommand*\linenomathpatch{\@ifstar{\linenomathpatch@AMS}{\linenomathpatch@}}
  \newcommand*\linenomathpatch@[1]{
    \expandafter\pretocmd\csname #1\endcsname {\linenomathWithnumbers}{}{}
    \expandafter\pretocmd\csname #1*\endcsname{\linenomathWithnumbers}{}{}
    \expandafter\apptocmd\csname end#1\endcsname {\endlinenomath}{}{}
    \expandafter\apptocmd\csname end#1*\endcsname{\endlinenomath}{}{}
  }
  \newcommand*\linenomathpatch@AMS[1]{
    \expandafter\pretocmd\csname #1\endcsname {\linenomathWithnumbersAMS}{}{}
    \expandafter\pretocmd\csname #1*\endcsname{\linenomathWithnumbersAMS}{}{}
    \expandafter\apptocmd\csname end#1\endcsname {\endlinenomath}{}{}
    \expandafter\apptocmd\csname end#1*\endcsname{\endlinenomath}{}{}
  }
  \let\linenomathWithnumbersAMS\linenomathWithnumbers
  \patchcmd\linenomathWithnumbersAMS{\advance\postdisplaypenalty\linenopenalty}{}{}{}
  \linenomathpatch{equation}
  \linenomathpatch*{gather}
  \linenomathpatch*{multline}
  \linenomathpatch*{align}
  \linenomathpatch*{alignat}
  \linenomathpatch*{flalign}
\fi
\makeatother
 
\ifpdf
\hypersetup{
  pdftitle={Efficient solution of parameter identification problems with
            \texorpdfstring{$H^1$}{H\textasciicircum1} regularization},
  pdfauthor={Jan Blechta, Oliver G. Ernst},
  pdfsubject={MSC2020 65F08, 65F22, 65N21},
  pdfkeywords={inverse problem,
               parameter identification,
               \texorpdfstring{$H^1$}{H\textasciicircum1} regularization,
               preconditioning,
               electrical resistivity tomography},
}
\fi
\title{%
  Efficient solution of parameter identification problems with $H^1$ regularization\thanks{%
    Submitted September~6, 2022, revised August~22, 2023.
    \funding{This work was supported by the European Union (EU) -- European Social Fund (ESF)
             and the Free State of Saxony, project GEOSax, grant 100310486.}}
}
\author{%
    Jan Blechta\thanks{%
    Faculty of Mathematics and Physics, Charles University, 186\;75 Prague, Czech Republic
    (\email{blechta@karlin.mff.cuni.cz}).}
    \and
    Oliver G.\ Ernst\thanks{%
    Department of Mathematics, Chemnitz University of Technology, Chemnitz, 09126 Germany
    (\email{oernst@math.tu-chemnitz.de}).}
}
\headers{Efficient parameter identification with $H^1$ regularization}
        {J.~Blechta and O.\ G.~Ernst}

\begin{document}
\maketitle

\begin{abstract}
  We consider the identification of spatially distributed parameters under $H^1$~regularization.
  Solving the associated minimization problem by Gauss--Newton iteration results in linearized problems to be
  solved in each step that can be cast as boundary value problems involving a~low-rank modification of the
  Laplacian.
  Using algebraic multigrid as a~fast Laplace solver, the Sherman--Morrison--Woodbury formula can be employed
  to construct a~preconditioner for these linear problems which exhibits excellent scaling w.r.t.\ the
  relevant problem parameters.
  We first develop this approach in the functional setting, thus obtaining a~consistent methodology for selecting
  boundary conditions that arise from the $H^1$ regularization.
  We then construct a~method for solving the discrete linear systems based on combining any fast Poisson solver
  with the Woodbury formula.
  The efficacy of this method is then demonstrated with scaling experiments.
  These are carried out for a~common nonlinear parameter identification problem arising in electrical resistivity
  tomography.
\end{abstract}



\tableofcontents

\section{Introduction}

The problem of reconstructing a~distributed parameter by the standard output least squares approach leads, after discretization, to a~system of (nonlinear) algebraic equations which is typically solved using Newton-type methods, requiring the solution of a~linearized problem at each step.
For Gauss--Newton iteration, this linearized problem involves the Jacobian of the parameter-to-observation map, resulting in a~discrete least squares problem.
In the common setting where a~high-dimensional unknown parameter is to be reconstructed from a~small number of observations, this least squares problem is typically rank-deficient.
The underdetermined nature combined with the ill-posedness of the underlying continuous inverse problem make it necessary to regularize the least squares formulation by adding a~penalty term, usually involving norms of spatial derivatives of the unknown.
In the high-dimensional setting the linearized problems are also solved iteratively, usually by some variant of Krylov subspace projection methods adapted to least squares problems.
The ill-posed nature of the underlying inverse problem as well as the spectral distribution of the regularization operator combine to make the preconditioning of the least squares iteration highly challenging and many strategies have been proposed.
For parameter identification problems arising from partial differential equations, the Jacobian is typically a~compact operator (cf.\ \cite{Worthen2012}, \cite[Theorem 4.21]{ColtonKress2013}), and it is known that Krylov subspace methods such as LSQR converge very slowly for the discretized problem.
Tikhonov regularization by adding an $L^2$-norm penalty term changes the setting to a~compact perturbation of the identity, which in turn is fundamentally favorable for fast convergence of Krylov subspace iterations
\cite{Winther1980,moret1997,EiermannErnst2001,blechta2021}.
The spectral properties of the iteration matrix, however, become much more challenging when regularization involving smoothing terms are employed, leading to large Gauss--Newton inner iteration counts.

In this work we address the efficient solution of the nonlinear least squares problem arising from distributed parameter estimation problems regularized by the $H^1$ norm, sometimes referred to as  \emph{smoothness regularization}.
The function to be reconstructed from observations is represented as piecewise constant w.r.t.\ a~given triangulation of the domain and, following an idea proposed by~Schwarzbach and Haber~\cite{schwarzbach-haber-2013}, define its gradient in the regularization term by duality.
We derive this formulation in the continuous setting by expressing the Gauss--Newton updates as the solution of boundary value problems, which are then discretized using a~mixed finite element method.
After discretization, the linear systems arising in each Gauss--Newton step have a~saddle-point structure and are solved by  preconditioned MINRES \cite{PaigeSaunders1975,Fischer1996} iteration.
Our preconditioning strategy is based on a~known spectrally equivalent preconditioner for the Laplacian, which enters the problem by way of the regularization term, combined with an application of the Sherman--Morrison--Woodbury formula to account for the low-rank perturbation arising from the data misfit term.
As a~result, we obtain a~mesh-independent solver for the nonlinear least squares problem which is also robust w.r.t.\ a~large range of regularization parameters.

Background references for PDE-based nonlinear parameter identification problems are \cite{Vogel2002} and, with a~focus on geoelectromagnetic exploration problems, \cite{Haber2014}.
While less common for parameter identification, $H^1$ regularization is also used in optimal control problems involving control and state constraints \cite{BarkerEtAl2016,ChenHuang2017}.
The extensive literature on Krylov projection methods for least squares problems is summarized in \cite{Bjorck1996,Bjorck2015}.
Finite-precision effects are analyzed in \cite{BjorckEtAL1998} and more recent developments include extensions of these methods to the class of \emph{symmetric quasi-definite} problems in \cite{OrbanArioli2017}, an analysis of LSQR for compact operators in Hilbert space in \cite{CarusoNovati2019}, and an extensive numerical comparison of the state of the art in Krylov methods and preconditioners for sparse linear least squares problems in \cite{GouldScott2017}.
A~popular construction principle for preconditioning matrices of saddle-point structure is based on the observation that suitable block triangular and block diagonal preconditioners result in a~system matrix with a~minimal polynomial of degree two or three \cite{KlawonnStarke1999,murphy-golub-wathen-2000,ipsen-2001}, for which Krylov subspace projection will return the exact solution in the same number of steps.
A~more comprehensive review of operator preconditioning techniques with special emphasis on mixed discretizations and saddle point problems can be found in \cite{MardalWinther2011}.
A~large class of preconditioning techniques for general least squares problems are based on incomplete factorizations \cite{BenziTuma2003,BruEtAl2014,ArioliDuff2015,ScottTuma2016,Scott2017,CerdanEtAl2020} as well as inner-outer iteration \cite{MorikuniHayami2013}.
Closer to the approach proposed in this work, the idea of using a~suitable Laplace preconditioner for variational inverse problems involving a~compact operator, when the Laplacian is used as a~regularization for the normal equations, is explored in~\cite{HankeVogel1999} (cf.\ also~\cite{HankeVogel1998}).
A~refinement of this approach is described in~\cite{JacobsenEtAl2003} and a~further variant proposed in~\cite{Bunse-GerstnerEtAl2006}.
Image restoration problems are also close to our setting in that the origin of the least squares problem is a~continuous inverse problem and regularization is a~necessity.
However, given that there the unknown is an unblurred image, the basic problem is typically not one of least squares, since there are typically as many measurements (pixel values) as unknowns.
Once regularization by penalty terms is added, however, the formulation is typically that of a~minimization problem \cite{BerishaNagy2014,ChungGazzola2021TR}, and a~successful approach here is the class of \emph{hybrid projection methods} \cite{ChungPalmer2015,ChungEtAl2015}.
Particularly in connection with statistical inverse problems, using preconditioners derived from covariance matrices have recently drawn increased attention \cite{CalvettiSomersalo2005,CalvettiEtAl2018}.
Finally, methods employing the Sherman--Morrison--Woodbury formula for constructing preconditioners have been considered by Yin~\cite{Yin2009} and Benzi and Faccio~\cite{benzi-faccio-2022}.
In~\cite{Yin2009}, a~recursive factorization technique is employed to apply a~preconditioner for Tikhonov-regularized least squares problems with a~Euclidean penalty term.
Benzi and Faccio~\cite{benzi-faccio-2022} discuss preconditioning strategies for linear systems with a~matrix of the form $\boldsymbol{A} + \gamma\boldsymbol{U}\boldsymbol{U}^\top$ with a~tall-and-skinny matrix~$\boldsymbol{U}$,
a~class of problems which includes that addressed in this paper.
It is reported there that efforts to construct a~preconditioner using the Woodbury formula with an approximation of the factor $\boldsymbol{A}^{-1}$ occurring therein proved unsuccessful~\cite[p.~4]{benzi-faccio-2022}.
By contrast, we will demonstrate this approach for our problem, in which $\boldsymbol{A}$ is a~discrete
 Laplacian, to be quite effective.
It is fair to remark that~\cite{benzi-faccio-2022} considered a~broader class of problems,
including the case of singular~$\boldsymbol{A}$.

The structure of the paper is as follows: \Cref{sec:problem} introduces the problem setting of $H^1$-regularized parameter estimation, derives the operator equations to be solved in each Gauss--Newton step, and briefly presents its mixed discretization based on an inf-sup stable mixed discretization for the Poisson equation.
\Cref{sec:la} presents three variants of the solution algorithm:
(i) direct approach based on the use of the Woodbury formula and a~factorization of the Laplacian,
(ii) MINRES iteration preconditioned by the Woodbury formula and a~suitable Laplace preconditioner,
(iii) a~simplified variant of~(ii) that omits the low-rank modification in the Woodbury formula.
\Cref{sec:experiments} contains an~extensive numerical illustration, in which
our solution approach is applied to an electrical resistivity tomography
problem from geophysical exploration in two and three space dimensions.
Realistic measurement setups are considered involving up to thousands of observational data points.
The algorithm is seen to perform efficiently and robustly across a~variety of settings.
Finally, in \cref{sec:outro} we summarize our findings and indicate further aspects to be investigated
in subsequent research.

\section{Problem formulation}\label{sec:problem}
We consider the output-least-squares formulation for estimating a~distributed parameter $m \in L^2(\Omega)$ defined on a~bounded domain $\Omega\subset{\mathbb R}^d$ and a~(typically nonlinear) parameter-to-observation map $g\colon L^2(\Omega)\to{\mathbb R}^M$ assigning to each parameter $m$ a~set of~$M$ observations, from which $m$ is to be reconstructed by minimizing the misfit $\sum_{i=1}^M \snorm{g^i(m) - {g^i_{\mathrm{obs}}}}^2$ w.r.t.\ a~vector ${\bm{g}_{\mathrm{obs}}}=\{{g^i_{\mathrm{obs}}}\}_{i=1}^M$ of observations.
For example, when $\exp m>0$ is the diffusion coefficient of an elliptic forward problem, this minimization is an ill-posed and severely underdetermined problem, which can be addressed by adding a~regularizing penalty term to the data misfit functional.
In this work we develop efficient computational methods for determining $m$ when the regularization term is the $H^1$ norm, a~common device for promoting smoothness of the reconstructed function.
This leads to the task of minimizing the objective function
\begin{equation} \label{eq:ls0}
  \sum_{i=1}^M \snorm{g^i(m) - {g^i_{\mathrm{obs}}}}^2 + \beta \int_\Omega \snorm{\nabla (m - {m_\mathrm{ref}})}^2,
\end{equation}
where $\beta>0$ is a~regularization parameter and ${m_\mathrm{ref}}$ denotes a~reference or background value for the unknown $m$.
The regularization thus penalizes the gradient of the deviation from the known background value ${m_\mathrm{ref}}$, a~common setting in, e.g., geophysical inverse problems.
Consequently, in order for the regularization term to make sense, this formulation, which we shall weaken in the following, would require $m-{m_\mathrm{ref}}$ to lie in the smaller space $H^1(\Omega)\subset L^2(\Omega)$.

To develop a~Gauss--Newton iteration for the minimization of~\cref{eq:ls0}, we will reformulate the first-order optimality condition as a~set of normal equations in the function space setting.
The gradient acting on $m$ then becomes a~Laplacian, for which spectrally equivalent preconditioners are available, allowing efficient iterative solution of the linearized equation in each Gauss--Newton step.
In addition, we recast the optimality equations in a~mixed formulation, which is well-defined also for $m \in L^2(\Omega)$, where the gradient in the regularization term is defined by duality.

\subsection{Assumptions and notation}
We assume that $\Omega\subset{\mathbb R}^d$ is a~bounded Lipschitz domain with boundary partitioned into ${\partial\Omega} = \overline{\Gamma_{\mathrm D}} \cup \overline{\Gamma_{\mathrm N}}$, with ${\Gamma_{\mathrm N}}$, ${\Gamma_{\mathrm D}}$ open and disjoint. For simplicity we assume $\snorm{{\Gamma_{\mathrm D}}}>0$ to exclude the pure Neumann problem.
We denote by $L^2(\Omega)$ the space of measurable functions $f\colon\Omega\to{\mathbb R}$ with finite norm
\begin{align*}
  \norm{f}_2 \coloneqq \biggl( \int_\Omega \snorm{f}^2 \biggr)^{\frac12}
\end{align*}
and by~$H^1(\Omega)$ the Sobolev space of such functions $f\in L^2(\Omega)$ with finite norm
\begin{align*}
  \norm{f}_{1,2} \coloneqq \biggl( \int_\Omega \snorm{\nabla f}^2 + \snorm{f}^2 \biggr)^{\frac12}.
\end{align*}
The subspace of functions vanishing on~${\Gamma_{\mathrm D}}$ is denoted by $H^1_{\Gamma_{\mathrm D}}(\Omega) \subset H^1(\Omega)$ and
\begin{align*}
  \norm{\nabla f}_2 \coloneqq \biggl( \int_\Omega \snorm{\nabla f}^2 \biggr)^{\frac12}
\end{align*}
is a~norm on~$H^1_{\Gamma_{\mathrm D}}(\Omega)$ that is equivalent to~$\norm{{\mkern 2mu\cdot\mkern 2mu}}_{1,2}$.
The space $H(\operatorname{div};\Omega)$ consists of all vector fields ${\vec{\mathsf{f}}}\colon\Omega\to{\mathbb R}^d$
such that $\snorm{{\vec{\mathsf{f}}}}\in L^2(\Omega)$ and $\operatorname{div}{\vec{\mathsf{f}}}\in L^2(\Omega)$ and it is equipped
with the norm
$
  \norm{{\vec{\mathsf{f}}}}_{\operatorname{div}{}}
  \coloneqq
  \norm{\snorm{{\vec{\mathsf{f}}}}}_2 + \norm{\operatorname{div}{\vec{\mathsf{f}}}}_2
$.
The subspace of~$H(\operatorname{div};\Omega)$
consisting of vector fields with vanishing normal trace on~${\Gamma_{\mathrm N}}$ is denoted
by~$H_{{\Gamma_{\mathrm N}}}(\operatorname{div};\Omega)$.

Next we assume that the parameter-to-observation map~$g$ is given
and G\^{a}teaux-differentiable with
the derivative denoted by $J\colon L^2(\Omega) \to [L^2(\Omega)']^M$ so that
\begin{align*}
  \langle J(m), \delta m \rangle = \biggl[ \frac{\mathrm{d}}{\mathrm{d}t}{g(m+t\delta m)} \biggr] _{\mkern 1mu \vrule height 2ex\mkern2mu t=0}
  \qquad
  m,\delta m\in L^2(\Omega).
\end{align*}
Therefore the mapping $\delta m \mapsto \langle J(m), \delta m \rangle$ is assumed
to be linear and bounded.
The individual components of $g$ and $J$ are denoted by $g^i\colon L^2(\Omega)\to{\mathbb R}$ and $J^i\colon L^2(\Omega)\to L^2(\Omega)'$, respectively, so that
\begin{align*}
  \langle J^i(m),\delta m \rangle = \biggl[ \frac{\mathrm{d}}{\mathrm{d}t}{g^i(m+t\delta m)} \biggr] _{\mkern 1mu \vrule height 2ex\mkern2mu t=0}
  \qquad
  m,\delta m\in L^2(\Omega),
  \quad i=1,2,\ldots,M.
\end{align*}

\subsection{Primal and mixed regularized least squares formulation}
To simplify the following expressions, we rescale the regularized least squares functional~\cref{eq:ls0} by $1/\beta>0$ and obtain the objective function
\begin{align} \label{eq:ls}
  \Phi_\beta(m)
  =
  \tfrac1\beta \sum_{i=1}^M \snorm{g^i(m)-{g^i_{\mathrm{obs}}}}^2 + \int_\Omega \snorm{\nabla(m-{m_\mathrm{ref}})}^2,
  \qquad m-{m_\mathrm{ref}}\in H^1_{\Gamma_{\mathrm D}}(\Omega).
\end{align}
Besides requiring the deviation $m-{m_\mathrm{ref}}$ to lie in the smoother space $H^1(\Omega)$, we impose an essential boundary condition on the portion ${\Gamma_{\mathrm D}}$ of the boundary of the domain $\Omega$.
This is a~modeling decision, which depends on the type of assumptions or a~priori information available on the unknown parameter $m$; in this case $m$ is assumed to coincide with the background value ${m_\mathrm{ref}}$ on ${\Gamma_{\mathrm D}}$.
As we will see below, this choice also implicitly imposes a~natural boundary condition on ${\Gamma_{\mathrm N}}$.

Taking the first variation (G\^{a}teaux derivative) of~\cref{eq:ls} in a~direction $\phi \in H^1_{\Gamma_{\mathrm D}}(\Omega)$ and setting it to zero, we arrive at the first-order necessary optimality condition for minimizing \cref{eq:ls}:
\begin{equation}  \label{eq:el}
  \begin{multlined}
    \text{Find } m\in H^1(\Omega) \text{ such that } m-{m_\mathrm{ref}}\in H^1_{\Gamma_{\mathrm D}}(\Omega) \text{ and}
    \\
    \hspace{2em}
    \tfrac1\beta
    \sum_{i=1}^M (g^i(m)-{g^i_{\mathrm{obs}}})\, \langle J^i(m),\phi \rangle
    + \int_\Omega \nabla(m-{m_\mathrm{ref}})\cdot\nabla\phi
    = 0
    \hspace{2em}
    \\
    \text{for all $\phi\in H^1_{\Gamma_{\mathrm D}}(\Omega)$.}
  \end{multlined}
\end{equation}
Assuming sufficient regularity, \cref{eq:el} can be interpreted as a~weak formulation of the boundary value problem
\begin{subequations} \label{eq:strong}
\begin{alignat}{2}
  \label{eq:stronga}
  \smash[b]{
    \tfrac1\beta
    \sum_{i=1}^M
    (g^i(m)-{g^i_{\mathrm{obs}}}) \, (J^i(m))'
    -\Delta(m-{m_\mathrm{ref}})
  }
  &= 0
  &\qquad&\text{ in } \Omega,
  \\[.2em]
  m-{m_\mathrm{ref}} &= 0
  &\qquad&\text{ on } {\Gamma_{\mathrm D}},
  \\
  \tfrac{\partial}{\partial\mathrm{n}}(m-{m_\mathrm{ref}}) &= 0
  &\qquad&\text{ on } {\Gamma_{\mathrm N}}.
\end{alignat}
\end{subequations}
Here $(J^i(m))'\in L^2(\Omega)$ denotes the Riesz representer of $J^i(m)\in L^2(\Omega)'$, i.e.,
\begin{align}
  \label{eq:riesz}
  \langle J^i(m), \phi \rangle = \int_\Omega (J^i(m))' \phi
  \qquad \text{for all }\phi \in L^2(\Omega).
\end{align}
Recall that $g^i(m)-{g^i_{\mathrm{obs}}}$ is a~number for any fixed $m$.
Hence, for a~fixed $m$, the first term
in~\cref{eq:stronga} is an $L^2(\Omega)$-function in the present setting.

Gauss--Newton linearization of~\cref{eq:el} is obtained by applying Newton's method to~\cref{eq:el} and neglecting the Hessian of~$g$,
which is given by
\begin{equation*}
  \label{eq:hessian}
  \langle H(m) \phi, \delta m \rangle =
  \biggl[ \frac{\mathrm{d}}{\mathrm{d}t}{\langle J(m+t\delta m), \phi \rangle} \biggr] _{\mkern 1mu \vrule height 2ex\mkern2mu t=0}
  \qquad
  m\in H^1(\Omega),\quad \phi,\delta m\in H^1_{\Gamma_{\mathrm D}}(\Omega).
\end{equation*}
Given an~initial value~$m$, associated model-generated responses~$g_m^i = g^i(m)$,
and the derivatives~$J_{m}^i = J^i(m)$,\; $i=1,2,\ldots,M$,
one step of Gauss--Newton iteration determines an update $m+\delta m$ by solving
the following problem:
\begin{equation} \label{eq:gn}
\begin{multlined}
  \text{Find } \delta m\in H^1(\Omega) \text{ such that } \delta m+m-{m_\mathrm{ref}}\in H^1_{\Gamma_{\mathrm D}}(\Omega) \text{ and }
  \\
  \hspace{2em}
  \begin{aligned}
    &
    \tfrac1\beta
    \sum_{i=1}^M \langle J_{m}^i,\delta m\rangle \, \langle J_{m}^i,\phi\rangle
    + \int_\Omega \nabla\delta m\cdot\nabla\phi
    \\
    &\hspace{2em}
    =
    -\tfrac1\beta
    \sum_{i=1}^M (g_m^i-{g^i_{\mathrm{obs}}}) \, \langle J_{m}^i,\phi\rangle
    -\int_\Omega \nabla(m-{m_\mathrm{ref}})\cdot\nabla\phi
  \end{aligned}
  \hspace{2em}
  \\
  \text{for all } \phi\in H^1_{\Gamma_{\mathrm D}}(\Omega).
\end{multlined}
\end{equation}
The variational equation \cref{eq:gn} in turn is a~weak formulation of the
boundary value problem for the Gauss--Newton correction $\delta m$
\begin{subequations} \label{eq:gn_strong}
\begin{alignat}{2}
  \label{eq:gn_strong_a}
  \Biggl[
    \tfrac1\beta\sum_{i=1}^M {J_{m}^i}'\, \langle J_{m}^i, {\mkern 2mu\cdot\mkern 2mu} \rangle
    -\Delta
  \Biggr] \delta m
  &=
  -\tfrac1\beta\sum_{i=1}^M {J_{m}^i}'\, (g_m^i-{g^i_{\mathrm{obs}}})
  +\Delta(m-{m_\mathrm{ref}})
  &\enspace&\text{in } \Omega,
  \\
  \delta m &= -(m-{m_\mathrm{ref}})
  &\enspace&\text{on } {\Gamma_{\mathrm D}},
  \\[.5em]
  \tfrac{\partial}{\partial\mathrm{n}}\delta m &= -\tfrac{\partial}{\partial\mathrm{n}}(m-{m_\mathrm{ref}})
  &\enspace&\text{on } {\Gamma_{\mathrm N}}.
\end{alignat} \end{subequations}
In view of \cref{eq:riesz}, the first operator in~\cref{eq:gn_strong_a} acting on $\delta m$
can be expressed as
\begin{align*}
  \Biggl[
    \sum_{i=1}^M {J_{m}^i}'\, \langle J_{m}^i, \delta m \rangle
  \Biggr] (x)
  =
  \int_\Omega
    \Biggl[
      \sum_{i=1}^M {J_{m}^i}'(x)\, {J_{m}^i}'(y)
    \Biggr]
    \delta m(y)\, \mathrm{d}y,
\end{align*}
i.e., as a~finite-rank integral operator with kernel
$\sum_{i=1}^M {J_{m}^i}'(x)\, {J_{m}^i}'(y) \in L^2(\Omega\times\Omega)$.

To weaken the regularity requirements on $m$ we next recast problem~\cref{eq:gn} in a~mixed formulation by
introducing the flux variable ${\vec{\zeta}}\coloneqq \nabla(\delta m+m-{m_\mathrm{ref}})$:
\begin{equation} \label{eq:mixed}
  \begin{multlined}
    \text{Find } ({\vec{\zeta}},\delta m)\in H_{{\Gamma_{\mathrm N}}}(\operatorname{div};\Omega)\times L^2(\Omega) \text{ such that}
    \\[.3em]
    \hspace{2em}
    \begin{aligned}
       \int_\Omega {\vec{\zeta}} \cdot {\vec{\psi}} + \int_\Omega \delta m\, \operatorname{div} {\vec{\psi}}
      &=
      -\int_\Omega (m-{m_\mathrm{ref}})\, \operatorname{div}{\vec{\psi}},
      \\
      \tfrac1\beta
      \sum_{i=1}^M \langle J_{m}^i,\delta m\rangle \, \langle J_{m}^i,\phi\rangle
      - \int_\Omega \phi\, \operatorname{div}{\vec{\zeta}}
      &=
      -\tfrac1\beta
      \sum_{i=1}^M (g_m^i-{g^i_{\mathrm{obs}}}) \, \langle J_{m}^i,\phi\rangle
    \end{aligned}
    \hspace{2em}
    \\[.1em]
    \text{for all } {\vec{\psi}}\in H_{{\Gamma_{\mathrm N}}}(\operatorname{div};\Omega) \text{ and } \phi\in L^2(\Omega).
  \end{multlined}
\end{equation}
Introducing the operators
\begin{subequations} \label{eq:op}
\begin{alignat}{2} \label{eq:opM}
  \langle Q{\vec{\zeta}},{\vec{\psi}} \rangle
  &\coloneqq \int_\Omega {\vec{\zeta}}\cdot{\vec{\psi}}
  &\qquad& {\vec{\zeta}},{\vec{\psi}}\in H_{{\Gamma_{\mathrm N}}}(\operatorname{div};\Omega), \\
  \label{eq:opD}
  \langle D{\vec{\zeta}},\phi \rangle
  &\coloneqq \int_\Omega \phi\,  \operatorname{div}{\vec{\zeta}}
  &\qquad& {\vec{\zeta}}\in H_{{\Gamma_{\mathrm N}}}(\operatorname{div};\Omega),\, \phi\in L^2(\Omega),
\end{alignat}
\end{subequations}
we can rewrite~\cref{eq:mixed} in the block operator form
\begin{equation} \label{eq:block_mixed}
  \begin{bmatrix}
      Q & D' \\
            D & -\frac1\beta J_{m}' J_{m}
  \end{bmatrix}
  \begin{bmatrix}
    {\vec{\zeta}} \\
    \delta m
  \end{bmatrix}
  =
  \begin{bmatrix}
    -D'(m-{m_\mathrm{ref}}) \\
     \frac1\beta J_{m}' (\bm{g}_{m}-{\bm{g}_{\mathrm{obs}}})
  \end{bmatrix},
\end{equation}
where the occurrences of $J_{m}$ and $J_{m}'$ are expressed using duality as
\begin{alignat*}{2}
  \langle J_{m}' J_{m} \delta m,\phi\rangle
  &= \sum_{i=1}^M
  \langle J_{m}^i,\delta m\rangle \, \langle J_{m}^i,\phi\rangle
  &\qquad& \delta m,\phi\in L^2(\Omega),
  \\
  \langle J_{m}' (\bm{g}_{m}-{\bm{g}_{\mathrm{obs}}}),\phi\rangle
  &= \sum_{i=1}^M
  (g_m^i-{g^i_{\mathrm{obs}}}) \, \langle J_{m}^i,\phi\rangle
  &\qquad& \phi\in L^2(\Omega).
\end{alignat*}
In an analogous way, by defining the operator
\begin{alignat*}{2}
  \langle L\delta m,\phi \rangle
  &\coloneqq \int_\Omega \nabla\delta m\cdot\nabla\phi
  &\qquad& \delta m,\phi\in H^1_{\Gamma_{\mathrm D}}(\Omega),
\end{alignat*}
we can rewrite the primal formulation~\cref{eq:gn} as the operator equation
\begin{equation} \label{eq:block_primal}
  \bigl( L+\tfrac1\beta J_{m}' J_{m} \bigr) \delta m
  = -L(m-{m_\mathrm{ref}}) -\tfrac1\beta J_{m}' (\bm{g}_{m}-{\bm{g}_{\mathrm{obs}}}).
\end{equation}
The primal formulation~\cref{eq:block_primal} can be seen as Schur complement reduction of the mixed formulation~\cref{eq:block_mixed}, in terms of which the Laplacian is represented as $L=D Q^{-1} D'$.
Indeed, block elimination of~${\vec{\zeta}}$ in~\cref{eq:block_mixed} gives
\begin{equation} \label{eq:block_schur}
  \bigl( D Q^{-1} D' + \tfrac1\beta J_{m}' J_{m} \bigr) \delta m
  = - D Q^{-1} D'(m-{m_\mathrm{ref}}) -\tfrac1\beta J_{m}' (\bm{g}_{m}-{\bm{g}_{\mathrm{obs}}}).
\end{equation}
We note that Schwarzbach and Haber~\cite[section~3.2.1]{schwarzbach-haber-2013} also formulated the $H^1$~regularization using the mixed formulation~\cref{eq:block_mixed}.
Their approach consisted of discretizing by lowest-order Raviart--Thomas elements and approximating $Q^{-1}$ by a~diagonal matrix in the Schur complement formulation~\cref{eq:block_schur}.
We will instead proceed by considering the mixed formulation~\cref{eq:block_mixed} and design solution strategy for this system.

\subsection{Finite element discretization}

The $H^1(\Omega)$-formulation~\cref{eq:block_primal} suggests an $H^1(\Omega)$-conforming discretization for the parameter~$m$, using, e.g., continuous Lagrange elements.
Instead, to allow for parameters $m \in L^2(\Omega)$ we will employ a~standard discretization of the mixed formulation~\cref{eq:mixed}.
Let us assume in the following that $\Omega$ is polyhedral so that we can consider its simplicial partitions $\mathcal{T}_h$.
Further let finite element spaces
$V_h\times Q_h \subset H_{{\Gamma_{\mathrm N}}}(\operatorname{div};\Omega)\times L^2(\Omega)$
be chosen as
\begin{equation*}
  V_h \times Q_h \coloneqq
  RT_k(\mathcal{T}_h,\, {\Gamma_{\mathrm N}}) \times \mathrm{d}\mkern-1mu P_k(\mathcal{T}_h)
\end{equation*}
for some order $k\in \mathbb N_0$, where $RT_k$ and $\mathrm{d}\mkern-1mu P_k$ denote the finite element spaces of Raviart--Thomas
and discontinuous Lagrange of order~$k$ counted such that $k=0$ corresponds to the lowest-order case.
This is an~inf-sup stable discretization for the Poisson equation in mixed formulation, i.e., the operator given by~\cref{eq:block_mixed} (or, equivalently, \cref{eq:mixed}) without the $J_{m}' J_{m}$ term; see~\cite{boffi-brezzi-fortin}.

Let $\{ {\vec{\psi}}^h_i \}_{i=1}^K$ and $\{ \phi^h_i \}_{i=1}^N$ denote bases of $V_h$ and $Q_h$, respectively, so that
\begin{alignat*}{2}
  \operatorname{span} \{ {\vec{\psi}}^h_i \}_{i=1}^K &= V_h,
  &\qquad
  K &= \dim V_h,
  \\
  \operatorname{span} \{ \phi^h_i \}_{i=1}^N &= Q_h,
  &\qquad
  N &= \dim Q_h.
\end{alignat*}
Inserting the basis elements into~\cref{eq:op} yields the matrices
\begin{subequations} \label{eq:oph}
\begin{alignat}{2}
  \label{eq:opMh}
  \boldsymbol{Q} &\in{\mathbb R}^{K\times K},
  &\qquad
  (\boldsymbol{Q})_{ij}
  &\coloneqq
  \langle Q{\vec{\psi}}^h_j,{\vec{\psi}}^h_i \rangle,
  \\
  \label{eq:opDh}
  \boldsymbol{D} &\in{\mathbb R}^{N\times K},
  &\qquad
  (\boldsymbol{D})_{ij}
  &\coloneqq
  \langle D{\vec{\psi}}^h_j,\phi^h_i \rangle.
\end{alignat}
\end{subequations}
Assuming $m,\,{m_\mathrm{ref}}\in Q_h$, these can then be expressed as
\begin{subequations} \label{eq:vectors}
\begin{alignat}{3}
  \label{eq:vectors_m}
  \bm{m} &\in{\mathbb R}^N,
  &\qquad
  m(x) &= \sum_{j=1}^N (\bm{m})_j \phi^h_j(x)
  &\quad x&\in\Omega,
  \\
  {\bm{m}_\mathrm{ref}} &\in{\mathbb R}^N,
  &\qquad
  {m_\mathrm{ref}}(x) &= \sum_{j=1}^N ({\bm{m}_\mathrm{ref}})_j \phi^h_j(x)
  &\quad x&\in\Omega.
  \intertext{
    We seek to determine
        $({\vec{\zeta}},\delta m) \in V_h\times Q_h$
    so that
  }
  \boldsymbol{\zeta} &\in{\mathbb R}^K,
  &\qquad
  {\vec{\zeta}}(x) &= \sum_{j=1}^K (\boldsymbol{\zeta})_j {\vec{\psi}}^h_j(x)
  &\quad x&\in\Omega,
  \\
  \label{eq:vectors_dm}
  \boldsymbol{\delta}\bm{m} &\in{\mathbb R}^N,
  &\qquad
  \delta m(x) &= \sum_{j=1}^N (\boldsymbol{\delta}\bm{m})_j \phi^h_j(x)
  &\quad x&\in\Omega.
\end{alignat}
Naturally, $g$ and $J$ are restricted to~$Q_h$, which gives rise
to the vector and the matrix
\begin{alignat}{2}
  \label{eq:opgh}
  \bm{g}_{\bm{m}} &\in{\mathbb R}^M,
  &\qquad
  (\bm{g}_{\bm{m}})_i &\coloneqq g^i(m),
  \\
  \label{eq:opJh}
  {\boldsymbol{J}}_{\bm{m}} &\in{\mathbb R}^{M\times N},
  &\qquad
  ({\boldsymbol{J}}_{\bm{m}})_{ij}
  &\coloneqq
  \langle J^i(m),\phi^h_j\rangle,
\end{alignat}
where $m\in Q_h$ on the right-hand sides is given by~\cref{eq:vectors_m}.
\end{subequations}
We thus arrive at the discrete counterpart of~\cref{eq:block_mixed}, the block linear system
\begin{equation} \label{eq:block_mixed_h}
  \begin{bmatrix}
     \boldsymbol{Q} & \boldsymbol{D}^\top \\
           \boldsymbol{D} & -\frac{1}{\beta} {\boldsymbol{J}}_{\bm{m}}^\top {\boldsymbol{J}}_{\bm{m}}
  \end{bmatrix}
  \begin{bmatrix}
    \boldsymbol{\zeta} \\
    \boldsymbol{\delta}\bm{m}
  \end{bmatrix}
  =
  \begin{bmatrix}
    -\boldsymbol{D}^\top(\bm{m}-{\bm{m}_\mathrm{ref}}) \\
     \tfrac{1}{\beta}{\boldsymbol{J}}_{\bm{m}}^\top(\bm{g}_{\bm{m}}-{\bm{g}_{\mathrm{obs}}})
  \end{bmatrix}.
\end{equation}

\section{Solution of the linear systems using the Woodbury formula}
\label{sec:la}

The linear system~\cref{eq:block_mixed_h} to be solved for the Gauss--Newton updates is a~low-rank perturbation of a~Poisson problem in the mixed formulation.
In this \namecref{sec:la} we employ the Woodbury matrix identity (see, e.g., \cite[section~2.1.4]{GolubVanLoan2013})
to construct algorithms for efficiently solving this system.
We first consider a~direct solution approach which can benefit from reusing the factorization for the unperturbed
problem.
As a~second approach, we propose two preconditioners for an iterative solution which can take advantage of any available efficient solution method for the unperturbed problem.

To this end we note that any formulation and discretization for the Laplacian which has an~efficient solution method can be used to proceed along the following lines.
Recall that the operators in the infinite-dimensional linear systems \cref{eq:block_mixed,,eq:block_primal,,eq:block_schur} are seen to be finite-rank perturbations of the Laplacian.
But in the following we will focus entirely on the discrete mixed formulation~\cref{eq:block_mixed_h}.

\subsection{Direct solution}

\afterpage{%
  \begin{algorithm2e}[b!]
\caption{Gauss--Newton with direct solver}
\label{alg:gn-direct}
\KwIn{\justifying%
      Parameter-to-observation map $\bm{g}\colon{\mathbb R}^N\to{\mathbb R}^M$,
      observational data ${\bm{g}_{\mathrm{obs}}}\in{\mathbb R}^M$,
      reference parameter ${\bm{m}_\mathrm{ref}} \in {\mathbb R}^N$,
      simplicial partition $\mathcal{T}_h$ of~$\Omega$,
      Neumann boundary ${\Gamma_{\mathrm N}}\subset{\partial\Omega}$,
      initial guess $\bm{m} \in {\mathbb R}^N$,
      regularization parameter $\beta > 0$}
\KwOut{Final $\bm{m} \in {\mathbb R}^N$}
Assemble matrices $\boldsymbol{Q} \in {\mathbb R}^{K\times K}$ and $\boldsymbol{D} \in {\mathbb R}^{N\times K}$
according to~\cref{eq:oph} \;
Use a~sparse direct solver
(e.g., sparse $\boldsymbol{L}\boldsymbol{D}\boldsymbol{L}^\top$)
to factorize
\nllabel{ln:ldlt}
\begin{align*}
  \boldsymbol{A} &=
  \begin{bmatrix}
     \boldsymbol{Q} & \boldsymbol{D}^\top \\
           \boldsymbol{D} & \boldsymbol{0}
  \end{bmatrix}
\end{align*}
\vskip-\belowdisplayskip\vskip\belowdisplayshortskip
\Repeat{happy}{
  Compute model response $\bm{g}_{\bm{m}} \in {\mathbb R}^M$
  and its derivative ${\boldsymbol{J}}_{\bm{m}} \in {\mathbb R}^{M\times N}$
  according to~\cref{eq:opgh,eq:opJh,eq:vectors_m}
  \nllabel{ln:observe}
  \;
  Use the factorization of $\boldsymbol{A}$ from above to construct $\boldsymbol{H}_{\bm{m}} \in {\mathbb R}^{N\times M}$ and $\boldsymbol{C}_{\beta,\bm{m}} \in {\mathbb R}^{M\times M}$
  such that
  \nllabel{ln:ldlt2}
  \begin{align*}
    \boldsymbol{H}_{\bm{m}} \coloneqq - \boldsymbol{P}_2 \boldsymbol{A}^{-1}
    \begin{bmatrix} \boldsymbol{0} \\ {\boldsymbol{J}}_{\bm{m}}^\top \end{bmatrix},
    \quad
    \boldsymbol{C}_{\beta,\bm{m}} \coloneqq \boldsymbol{I}_M + \tfrac{1}{\beta} {\boldsymbol{J}}_{\bm{m}} \boldsymbol{H}_{\bm{m}}
  \end{align*}
  \nl
  Solve the capacitance system directly,
  i.e., find $\bm{y} \in {\mathbb R}^{M}$ such that
  \nllabel{ln:capa}
  \begin{align*}
    \boldsymbol{C}_{\beta,\bm{m}} \bm{y} = {\boldsymbol{J}}_{\bm{m}} (\bm{m}-{\bm{m}_\mathrm{ref}}) - (\bm{g}_{\bm{m}}-{\bm{g}_{\mathrm{obs}}})
  \end{align*}
  \nl
  Compute $\boldsymbol{\delta}\bm{m} \in {\mathbb R}^N$ as
  \nllabel{ln:capa2}
  \begin{align*}
    \boldsymbol{\delta}\bm{m} \coloneqq
    -(\bm{m}-{\bm{m}_\mathrm{ref}})
    +\tfrac{1}{\beta} \boldsymbol{H}_{\bm{m}} \bm{y}
  \end{align*}
  \vskip-\belowdisplayskip\vskip\belowdisplayshortskip
  $\bm{m} \coloneqq \bm{m} + \boldsymbol{\delta}\bm{m}$ \;
}
   \end{algorithm2e}
}

We introduce the following matrices (cf.~\cref{eq:block_mixed_h}) which will be useful in constructing the solution schemes:
\begin{align}
  \label{eq:block_laplace}
   \boldsymbol{A} \coloneqq \begin{bmatrix}  \boldsymbol{Q} & \boldsymbol{D}^\top \\  \boldsymbol{D} & \boldsymbol{0} \end{bmatrix} \in {\mathbb R}^{(K+N)\times(K+N)},
   \qquad
   \begin{aligned}
     \boldsymbol{S}  &\coloneqq \boldsymbol{D}\boldsymbol{Q}^{-1}\boldsymbol{D}^\top \in {\mathbb R}^{N\times N},\\
     \boldsymbol{P}_2 &\coloneqq \begin{bmatrix} \boldsymbol{0} & \boldsymbol{I}_N \end{bmatrix} \in {\mathbb R}^{N\times(K+N)},
   \end{aligned}
\end{align}
where $\boldsymbol{I}_N$ denotes the $N \times N$ identity.
The solution of a~linear system with coefficient matrix $\boldsymbol{A}$ and right-hand side blocks $\bm{y}_1$ and $\bm{y}_2$ are related via the Schur complement~$\boldsymbol{S}$ as
\begin{align}
  \label{eq:schur}
  \boldsymbol{P}_2 \boldsymbol{A}^{-1} \begin{bmatrix} \bm{y}_1 \\ \bm{y}_2 \end{bmatrix}
  =
  \boldsymbol{S}^{-1} \bigl( \boldsymbol{D}\boldsymbol{Q}^{-1}\bm{y}_1 - \bm{y}_2 \bigr),
  \qquad \bm{y}_1 \in {\mathbb R}^{K},\; \bm{y}_2 \in {\mathbb R}^N.
\end{align}
The Schur complement matrix $\boldsymbol{S}$ is, in general, dense and hence linear systems with matrix~$\boldsymbol{S}$ are impractical to assemble and solve.
On the other hand, \cref{eq:schur}~implies that
\begin{align}
  \label{eq:schur2}
  \boldsymbol{S}^{-1} \bm{y}_2
  =
  - \boldsymbol{P}_2 \boldsymbol{A}^{-1} \begin{bmatrix} \bm{\mkern-1.5mu \mbox{\slshape 0}} \\ \bm{y}_2 \end{bmatrix},
  \qquad \bm{y}_2 \in {\mathbb R}^N,
\end{align}
i.e., the solution of the dense system $\boldsymbol{S}\bm{x}_2=-\bm{y}_2$ can be expressed as the
solution of the sparse saddle point system
\begin{align*}
  \SwapAboveDisplaySkip
  \boldsymbol{A} \begin{bmatrix} \bm{x}_1 \\ \bm{x}_2 \end{bmatrix}
  = \begin{bmatrix} \bm{\mkern-1.5mu \mbox{\slshape 0}} \\ \bm{y}_2 \end{bmatrix}.
\end{align*}
This is the setting we encounter in the Gauss--Newton update step, where
the second block of the solution of~\cref{eq:block_mixed_h} is needed and
the saddle point matrix~$\boldsymbol{A}_{\beta,\bm{m}}$ is a~low-rank
modification of~$\boldsymbol{A}$ in the second block:
\begin{align} \label{eq:block_sol}
  \boldsymbol{\delta}\bm{m} &= \boldsymbol{P}_2 \boldsymbol{A}_{\beta,\bm{m}}^{-1}
  \begin{bmatrix}
    -\boldsymbol{D}^\top(\bm{m}-{\bm{m}_\mathrm{ref}}) \\[.2em]
     \tfrac{1}{\beta}{\boldsymbol{J}}_{\bm{m}}^\top(\bm{g}_{\bm{m}}-{\bm{g}_{\mathrm{obs}}})
  \end{bmatrix},
  &
  \boldsymbol{A}_{\beta,\bm{m}} &\coloneqq \boldsymbol{A} - \frac{1}{\beta}
  \begin{bmatrix} \boldsymbol{0} \\ {\boldsymbol{J}}_{\bm{m}}^\top \end{bmatrix}
  \begin{bmatrix} \boldsymbol{0} & {\boldsymbol{J}}_{\bm{m}} \end{bmatrix}.
\end{align}
Using the Woodbury formula we may express $\boldsymbol{A}_{\beta,\bm{m}}^{-1}$ in terms of $\boldsymbol{A}^{-1}$ as
\begin{align*}
  \boldsymbol{A}_{\beta,\bm{m}}^{-1} &= \boldsymbol{A}^{-1} + \frac{1}{\beta}
  \boldsymbol{A}^{-1}
  \begin{bmatrix} \boldsymbol{0} \\ {\boldsymbol{J}}_{\bm{m}}^\top \end{bmatrix}
  \Biggl(
    \boldsymbol{I}_M -
    \frac{1}{\beta}
    \begin{bmatrix} \boldsymbol{0} & {\boldsymbol{J}}_{\bm{m}} \end{bmatrix}
    \boldsymbol{A}^{-1}
    \begin{bmatrix} \boldsymbol{0} \\ {\boldsymbol{J}}_{\bm{m}}^\top \end{bmatrix}
  \Biggr)^{-1}
  \begin{bmatrix} \boldsymbol{0} & {\boldsymbol{J}}_{\bm{m}} \end{bmatrix}
  \boldsymbol{A}^{-1},
\end{align*}
where $\boldsymbol{I}_M$ denotes the $M\times M$ identity.
Defining the matrix
\begin{gather} \label{eq:H}
  \boldsymbol{H}_{\bm{m}}
  \coloneqq
  \boldsymbol{S}^{-1} {\boldsymbol{J}}_{\bm{m}}^\top
  =
  - \boldsymbol{P}_2 \boldsymbol{A}^{-1} \begin{bmatrix} \boldsymbol{0} \\ {\boldsymbol{J}}_{\bm{m}}^\top \end{bmatrix}
  \in {\mathbb R}^{N\times M},
\intertext{and observing $\begin{bmatrix}\boldsymbol{0} & {\boldsymbol{J}}_{\bm{m}}\end{bmatrix}=  {\boldsymbol{J}}_{\bm{m}} \boldsymbol{P}_2$, we arrive at the expression for the matrix whose action is required in the update \cref{eq:block_sol}}
  \label{eq:woodbury2}
  \boldsymbol{P}_2 \boldsymbol{A}_{\beta,\bm{m}}^{-1}
  =
  \biggl(
    \boldsymbol{I}_N - \tfrac{1}{\beta}
    \boldsymbol{H}_{\bm{m}}
    \Bigl(
      \boldsymbol{I}_M +
      \tfrac{1}{\beta}
      {\boldsymbol{J}}_{\bm{m}}
      \boldsymbol{H}_{\bm{m}}
    \Bigr)^{-1}
    {\boldsymbol{J}}_{\bm{m}}
  \biggr)
  \boldsymbol{P}_2 \boldsymbol{A}^{-1}.
\end{gather}
Combining~\cref{eq:schur,eq:H}, we obtain for the unperturbed problem
\begin{align} \label{eq:woodbury_temp}
  \boldsymbol{P}_2 \boldsymbol{A}^{-1}
  \begin{bmatrix}
    -\boldsymbol{D}^\top(\bm{m}-{\bm{m}_\mathrm{ref}}) \\[.2em]  \tfrac{1}{\beta}{\boldsymbol{J}}_{\bm{m}}^\top(\bm{g}_{\bm{m}}-{\bm{g}_{\mathrm{obs}}})
  \end{bmatrix}
  &=
  -(\bm{m}-{\bm{m}_\mathrm{ref}}) - \tfrac{1}{\beta} \boldsymbol{H}_{\bm{m}} (\bm{g}_{\bm{m}}-{\bm{g}_{\mathrm{obs}}}).
\end{align}
\Cref{eq:block_sol,eq:woodbury2,eq:woodbury_temp} now yield an expression for the update vector as
\begin{align*}
  \boldsymbol{\delta}\bm{m} =
  -(\bm{m}-{\bm{m}_\mathrm{ref}})
  +\tfrac{1}{\beta} \boldsymbol{H}_{\bm{m}}
  \Bigl(
    \boldsymbol{I}_M +
    \tfrac{1}{\beta}
    {\boldsymbol{J}}_{\bm{m}}
    \boldsymbol{H}_{\bm{m}}
  \Bigr)^{-1}
  \Bigl(
    {\boldsymbol{J}}_{\bm{m}} (\bm{m}-{\bm{m}_\mathrm{ref}}) - (\bm{g}_{\bm{m}}-{\bm{g}_{\mathrm{obs}}})
  \Bigr).
\end{align*}
The computations for constructing this vector within a~complete Gauss--Newton minimization are summarized in \cref{alg:gn-direct}.
It requires a~single $\boldsymbol{L}\boldsymbol{D}\boldsymbol{L}^\top$ factorization of the large sparse matrix $\boldsymbol{A}$ (\cref{ln:ldlt}).
This is done once, prior to the nonlinear iteration, hence its computational cost is amortized over the nonlinear solution process.
On the other hand, the fill-in resulting in the factors of~$\boldsymbol{A}$, especially in~3D, makes application of~$\boldsymbol{A}^{-1}$ expensive with complexity considerably larger than $O(N)$.
This occurs $M$~times on \cref{ln:ldlt2}  and thus potentially becomes
a~bottleneck of the algorithm if $N$ and/or $M$ are large.
Once~$\boldsymbol{H}_{\bm{m}}$ is computed, the construction of the capacitance matrix $\boldsymbol{C}_{\beta,\bm{m}}$ on \cref{ln:ldlt2} can proceed very efficiently in~$O(M^2N)$ operations as a~BLAS Level~3 operation.
The dense solve on \cref{ln:capa} costs~$O(M^3)$ and can be efficiently performed by LAPACK.

The evaluation of the model response and its derivative on \cref{ln:observe} is assumed to be available as a~given function $\bm{m}\mapsto(\bm{g}_{\bm{m}},{\boldsymbol{J}}_{\bm{m}})$.
In many contexts, where the mapping is based on a~PDE model, the evaluation of~$\bm{g}_{\bm{m}}$ requires the solution
of a~\emph{forward PDE problem}, and the computation of the associated derivative~${\boldsymbol{J}}_{\bm{m}}$ can be performed efficiently using adjoint techniques.
This will be the case in the numerical examples presented in \cref{sec:experiments}.

\subsection{Iterative solution} \label{sec:iterative}

We consider the block-diagonal preconditioners
\begin{gather}
  \label{eq:pc-ideal}
  \boldsymbol{P}
  \coloneqq
  \begin{bmatrix}
    \boldsymbol{Q} &           \\
        & \boldsymbol{S}
  \end{bmatrix}
  \qquad
  \text{and}
  \qquad
  \boldsymbol{P}_{\beta,\bm{m}}
  \coloneqq
  \begin{bmatrix}
    \boldsymbol{Q} &                                 \\
        & \boldsymbol{S} + \tfrac1\beta{\boldsymbol{J}}_{\bm{m}}^\top{\boldsymbol{J}}_{\bm{m}}
  \end{bmatrix},
\end{gather}
with the Laplace Schur complement~$\boldsymbol{S}$ as in~\cref{eq:block_laplace}.
These are ``ideal'' preconditioners for~$\boldsymbol{A}$ and~$\boldsymbol{A}_{\beta,\bm{m}}$, respectively.
Indeed, the minimal polynomial of $\boldsymbol{A}\boldsymbol{P}^{-1}$ has degree at most~$3$ \cite[Proposition~1]{murphy-golub-wathen-2000} and, as a~consequence, minimum residual Krylov subspace iteration applied to $\boldsymbol{A}\boldsymbol{P}^{-1}$ converges in at most $3$ iterations, as shown by Murphy, Golub, and Wathen~\cite{murphy-golub-wathen-2000} (see also \cite[Theorem~2.2.3]{LiesenStrakos2013}).
This does not hold for $\boldsymbol{A}_{\beta,\bm{m}}\boldsymbol{P}_{\beta,\bm{m}}^{-1}$, but it is known that the spectrum of \smash{$\boldsymbol{A}_{\beta,\bm{m}}\boldsymbol{P}_{\beta,\bm{m}}^{-1}$}
is contained in $\smash{[-1,-\tfrac1\phi]\cup[1,\phi]}$, where $\phi=\smash{\frac{1+\sqrt{5}}{2}}$; see \cite[Theorem~4]{pearson-2013}.
This inclusion guarantees $2$-step linear convergence of MINRES
for $\boldsymbol{A}_{\beta,\bm{m}}\boldsymbol{P}_{\beta,\bm{m}}^{-1}$ independently of $M$, $N$, $\beta$, and the right-hand side;
see, e.g., \cite[section~3.1]{herzog-sachs-2015} or
\cite[section~4.2.4]{elman-silvester-wathen-2014}.

The action of~$\boldsymbol{P}^{-1}$ and $\boldsymbol{P}_{\beta,\bm{m}}^{-1}$ is essentially as expensive as that of~$\boldsymbol{A}^{-1}$ and $\boldsymbol{A}_{\beta,\bm{m}}^{-1}$, respectively, hence we seek a~good and inexpensive approximation of $\boldsymbol{P}^{-1}$ and $\boldsymbol{P}_{\beta,\bm{m}}^{-1}$.
Consider
\begin{gather*}
  \SwapAboveDisplaySkip
  \boldsymbol{\hat P}^{-1}
  \coloneqq
  \begin{bmatrix}
    \boldsymbol{\hat Q}^{-1} &           \\
              & \boldsymbol{\hat S}^{-1}
  \end{bmatrix}
  \qquad
  \text{and}
  \qquad
  \boldsymbol{\hat P}_{\beta,\bm{m}}^{-1}
  \coloneqq
  \begin{bmatrix}
    \boldsymbol{\hat Q}^{-1} &           \\
              & \boldsymbol{\hat S}_{\beta,\bm{m}}^{-1}
  \end{bmatrix},
  \shortintertext{where}
  \begin{aligned}
    \boldsymbol{\hat Q}^{-1}
    &\coloneqq
    (\operatorname{diag}\boldsymbol{Q})^{-1},
    \\
    \boldsymbol{\hat S}^{-1}
    &\coloneqq
    \texttt{AMG}(\boldsymbol{D}\boldsymbol{\hat Q}^{-1}\boldsymbol{D}^\top),
    \\
    \boldsymbol{\hat S}_{\beta,\bm{m}}^{-1}
    &\coloneqq
    \boldsymbol{\hat S}^{-1} - \tfrac1\beta\boldsymbol{\hat S}^{-1}{\boldsymbol{J}}_{\bm{m}}^\top
                (\boldsymbol{I}_M+\tfrac1\beta{\boldsymbol{J}}_{\bm{m}}\boldsymbol{\hat S}^{-1}{\boldsymbol{J}}_{\bm{m}}^\top)^{-1}
                {\boldsymbol{J}}_{\bm{m}}\boldsymbol{\hat S}^{-1}.
  \end{aligned}
\end{gather*}
The preconditioner $\boldsymbol{\hat P}^{-1}$ was introduced by Powell and Silvester~\cite{powell-silvester-2003}
for preconditioning the mixed Laplacian~$\boldsymbol{A}$.
We can employ this preconditioner also for~$\boldsymbol{A}_{\beta,\bm{m}}$, which is, in view of~\cref{eq:block_sol},
a~perturbation of~$\boldsymbol{A}$ by at most rank~$M$.
The expressions for $\boldsymbol{\hat P}_{\beta,\bm{m}}^{-1}$ follow easily by requiring, in analogy to \cref{eq:pc-ideal},
that $\boldsymbol{\hat S}_{\beta,\bm{m}}=\boldsymbol{\hat S}+\tfrac1\beta{\boldsymbol{J}}_{\bm{m}}^\top{\boldsymbol{J}}_{\bm{m}}$, and using the Woodbury matrix identity.

\Cref{alg:gn-iterative} summarizes the Gauss--Newton procedure based on
iterative solution of the linearized problems. It invokes
either \cref{alg:gnstep-iterative}, which employs \smash{$\boldsymbol{\hat P}_{\beta,\bm{m}}^{-1}$} as a~preconditioner,
or \cref{alg:gnstep-iterative-nw}, which uses $\boldsymbol{\hat P}^{-1}$.
The latter omits the correction due to the Woodbury formula,
hence bypasses the computations involving the capacitance matrix
$\boldsymbol{C}_{\beta,\bm{m}} = \boldsymbol{I}_M + \frac{1}{\beta} {\boldsymbol{J}}_{\bm{m}} \boldsymbol{\hat S}^{-1} {\boldsymbol{J}}_{\bm{m}}^\top$,
and thus results in a~less expensive preconditioner.
As it fails to account for the low-rank modification due to the data misfit term, it is
expected to deteriorate with increasing~$M$. We will confirm this
experimentally in \cref{sec:experiments}. Additionally we will see
that $\boldsymbol{\hat P}_{\beta,\bm{m}}^{-1}$, in contrast to~$\boldsymbol{\hat P}^{-1}$, provides robustness
w.r.t.~$\beta$; see \cref{fig:timing-beta}.
\begin{algorithm2e}[t]
\caption{Gauss--Newton with iterative solver}
\label{alg:gn-iterative}
\KwIn{\justifying%
      Parameter-to-observation map $\bm{g}\colon{\mathbb R}^N\to{\mathbb R}^M$,
      observational data ${\bm{g}_{\mathrm{obs}}}\in{\mathbb R}^M$,
      reference parameter ${\bm{m}_\mathrm{ref}} \in {\mathbb R}^N$,
      simplicial partition $\mathcal{T}_h$ of~$\Omega$,
      Neumann boundary ${\Gamma_{\mathrm N}}\subset{\partial\Omega}$,
      initial guess $\bm{m} \in {\mathbb R}^N$,
      regularization parameter $\beta > 0$}
\KwOut{Final $\bm{m} \in {\mathbb R}^N$}
Assemble matrices $\boldsymbol{Q} \in {\mathbb R}^{K\times K}$ and $\boldsymbol{D} \in {\mathbb R}^{N\times K}$
according to~\cref{eq:oph} \;
Prepare a~mass term preconditioner
\nonl
\begin{align*}
  \boldsymbol{\hat Q}^{-1} &\coloneqq (\operatorname{diag}\boldsymbol{Q})^{-1}
\end{align*}
\nl\nl
Prepare a~Schur complement preconditioner using an algebraic
blackbox, e.g., algebraic multigrid,
\begin{align*}
  \SwapAboveDisplaySkip
  \boldsymbol{\hat S}^{-1} &\coloneqq \texttt{AMG}(\boldsymbol{D}(\operatorname{diag}\boldsymbol{Q})^{-1}\boldsymbol{D}^\top)
\end{align*}
\vskip-\belowdisplayskip\vskip\belowdisplayshortskip
\Repeat{happy}{
  Compute model response $\bm{g}_{\bm{m}} \in {\mathbb R}^M$
  and its derivative ${\boldsymbol{J}}_{\bm{m}} \in {\mathbb R}^{M\times N}$
  according to~\cref{eq:opgh,eq:opJh,eq:vectors_m} \;
  Compute $\boldsymbol{\delta}\bm{m} \in {\mathbb R}^N$ using \cref{alg:gnstep-iterative}
  or \cref{alg:gnstep-iterative-nw} \;
  $\bm{m} \coloneqq \bm{m} + \boldsymbol{\delta}\bm{m}$ \;
}
 \end{algorithm2e}

\afterpage{%
  \begin{algorithm2e}[t]
\caption{MINRES with the Laplace--Woodbury preconditioner~$\boldsymbol{\hat P}_{\beta,\bm{m}}^{-1}$}
\label{alg:gnstep-iterative}
\KwIn{$\bm{m} \in {\mathbb R}^N$,
      $\bm{g}_{\bm{m}} \in {\mathbb R}^M$,
      ${\boldsymbol{J}}_{\bm{m}} \in {\mathbb R}^{M\times N}$,
      $\beta > 0$}
\KwOut{$\boldsymbol{\delta}\bm{m} \in {\mathbb R}^N$}
$\boldsymbol{\hat H}_{\bm{m}} \coloneqq \boldsymbol{\hat S}^{-1}{\boldsymbol{J}}_{\bm{m}}^\top$ \;
\nllabel{ln:H}
$\boldsymbol{C}_{\beta,\bm{m}} \coloneqq \boldsymbol{I}_M + \frac{1}{\beta} {\boldsymbol{J}}_{\bm{m}} \boldsymbol{\hat H}_{\bm{m}}$ \;
\nllabel{ln:C}
Compute Cholesky factor $\boldsymbol{L}_{\beta,\bm{m}} \in {\mathbb R}^{M\times M}$ such that
$\boldsymbol{L}_{\beta,\bm{m}}\boldsymbol{L}_{\beta,\bm{m}}^\top = \boldsymbol{C}_{\beta,\bm{m}}$ \;
\nllabel{ln:chol}
Run MINRES:
\nllabel{ln:minres}
\nonl
\begin{gather*}
  \SwapAboveDisplaySkip
  \begin{bmatrix} \boldsymbol{\zeta} \\ \boldsymbol{\delta}\bm{m} \end{bmatrix}
  \coloneqq
  \texttt{MINRES} \bigl( \boldsymbol{A}_{\beta,\bm{m}}, \boldsymbol{\hat P}_{\beta,\bm{m}}^{-1}, \bm{b}, \bm{x}_0 \bigr),
  \intertext{%
    where the system operator $\boldsymbol{A}_{\beta,\bm{m}}\colon{\mathbb R}^{K+N}\to{\mathbb R}^{K+N}$
    and the preconditioner $\boldsymbol{\hat P}_{\beta,\bm{m}}^{-1}\colon{\mathbb R}^{K+N}\to{\mathbb R}^{K+N}$
    are represented matrix-free by formulas
  }
  \begin{aligned}
    \boldsymbol{A}_{\beta,\bm{m}} \begin{bmatrix} \boldsymbol{\zeta} \\ \boldsymbol{\delta}\bm{m} \end{bmatrix}
    &=
    \begin{bmatrix}
       \boldsymbol{Q} \boldsymbol{\zeta} + \boldsymbol{D}^\top \boldsymbol{\delta}\bm{m}
      \\
      \boldsymbol{D} \boldsymbol{\zeta} - \frac{1}{\beta} ({\boldsymbol{J}}_{\bm{m}}^\top ({\boldsymbol{J}}_{\bm{m}}\boldsymbol{\delta}\bm{m}))
    \end{bmatrix},
    \\
    \boldsymbol{\hat P}_{\beta,\bm{m}}^{-1} \begin{bmatrix} \bm{y}_1 \\ \bm{y}_2 \end{bmatrix}
    &=
    \begin{bmatrix}
      \boldsymbol{\hat Q}^{-1} \bm{y}_1
      \\
      \boldsymbol{\hat S}^{-1} \bm{y}_2 - \frac{1}{\beta} (\boldsymbol{\hat H}_{\bm{m}}
        (\boldsymbol{L}_{\beta,\bm{m}}^{-\top} (\boldsymbol{L}_{\beta,\bm{m}}^{-1} (\boldsymbol{\hat H}_{\bm{m}}^\top \bm{y}_2))))
    \end{bmatrix}
  \end{aligned}
  \intertext{%
    and the right-hand side $\bm{b} \in {\mathbb R}^{K+N}$
    and the initial guess $\bm{x}_0 \in {\mathbb R}^{K+N}$
    are given by
  }
  \bm{b} \coloneqq
  \begin{bmatrix}
    -\boldsymbol{D}^\top (\bm{m}-{\bm{m}_\mathrm{ref}})
    \\
      \frac{1}{\beta} {\boldsymbol{J}}_{\bm{m}}^\top (\bm{g}_{\bm{m}}-{\bm{g}_{\mathrm{obs}}})
  \end{bmatrix},
  \quad
  \bm{x}_0 \coloneqq
  \begin{bmatrix}
    \bm{\mkern-1.5mu \mbox{\slshape 0}} \\ \bm{\mkern-1.5mu \mbox{\slshape 0}}
  \end{bmatrix}
\end{gather*}
   \end{algorithm2e}
  \begin{algorithm2e}[t]
\caption{MINRES with the Laplace preconditioner~$\boldsymbol{\hat P}^{-1}$}
\label{alg:gnstep-iterative-nw}
\KwIn{$\bm{m} \in {\mathbb R}^N$,
      $\bm{g}_{\bm{m}} \in {\mathbb R}^M$,
      ${\boldsymbol{J}}_{\bm{m}} \in {\mathbb R}^{M\times N}$,
      $\beta > 0$}
\KwOut{$\boldsymbol{\delta}\bm{m} \in {\mathbb R}^N$}
Run MINRES:
\nllabel{ln:minres-nw}
\nonl
\begin{gather*}
  \SwapAboveDisplaySkip
  \begin{bmatrix} \boldsymbol{\zeta} \\ \boldsymbol{\delta}\bm{m} \end{bmatrix}
  \coloneqq
  \texttt{MINRES} \bigl( \boldsymbol{A}_{\beta,\bm{m}}, \boldsymbol{\hat P}^{-1}, \bm{b}, \bm{x}_0 \bigr),
  \intertext{%
    where the system operator $\boldsymbol{A}_{\beta,\bm{m}}\colon{\mathbb R}^{K+N}\to{\mathbb R}^{K+N}$
    and the preconditioner $\boldsymbol{\hat P}_{\beta,\bm{m}}^{-1}\colon{\mathbb R}^{K+N}\to{\mathbb R}^{K+N}$
    are represented matrix-free by formulas
  }
  \begin{aligned}
    \boldsymbol{A}_{\beta,\bm{m}} \begin{bmatrix} \boldsymbol{\zeta} \\ \boldsymbol{\delta}\bm{m} \end{bmatrix}
    &=
    \begin{bmatrix}
       \boldsymbol{Q} \boldsymbol{\zeta} + \boldsymbol{D}^\top \boldsymbol{\delta}\bm{m}
      \\
      \boldsymbol{D} \boldsymbol{\zeta} - \frac{1}{\beta} ({\boldsymbol{J}}_{\bm{m}}^\top ({\boldsymbol{J}}_{\bm{m}}\boldsymbol{\delta}\bm{m}))
    \end{bmatrix},
    \\
    \boldsymbol{\hat P}^{-1} \begin{bmatrix} \bm{y}_1 \\ \bm{y}_2 \end{bmatrix}
    &=
    \begin{bmatrix}
      \boldsymbol{\hat Q}^{-1} \bm{y}_1
      \\
      \boldsymbol{\hat S}^{-1} \bm{y}_2
    \end{bmatrix}
  \end{aligned}
  \intertext{%
    and the right-hand side $\bm{b} \in {\mathbb R}^{K+N}$
    and the initial guess $\bm{x}_0 \in {\mathbb R}^{K+N}$
    are given by
  }
  \bm{b} \coloneqq
  \begin{bmatrix}
    -\boldsymbol{D}^\top (\bm{m}-{\bm{m}_\mathrm{ref}})
    \\
      \frac{1}{\beta} {\boldsymbol{J}}_{\bm{m}}^\top (\bm{g}_{\bm{m}}-{\bm{g}_{\mathrm{obs}}})
  \end{bmatrix},
  \quad
  \bm{x}_0 \coloneqq
  \begin{bmatrix}
    \bm{\mkern-1.5mu \mbox{\slshape 0}} \\ \bm{\mkern-1.5mu \mbox{\slshape 0}}
  \end{bmatrix}
\end{gather*}
   \end{algorithm2e}
}

The preconditioned MINRES procedures in \cref{alg:gnstep-iterative,alg:gnstep-iterative-nw}
correspond to different minimization problems
\begin{gather}
  \label{eq:minres}
  \doraisetag{-0.90cm}
  \begin{aligned}
  &\texttt{MINRES} \bigl( \boldsymbol{A}_{\beta,\bm{m}}, \boldsymbol{\hat P}_{\beta,\bm{m}}^{-1}, \bm{b}, \bm{x}_0 \bigr):
  &
  \norm{\bm{r}_k}_{\boldsymbol{\hat P}_{\beta,\bm{m}}^{-1}}
  &= \min_{p\in\mathcal{P}^0_k}
     \norm{p(\boldsymbol{A}_{\beta,\bm{m}}\boldsymbol{\hat P}_{\beta,\bm{m}}^{-1})\bm{r}_0}_{\boldsymbol{\hat P}_{\beta,\bm{m}}^{-1}},
  \\
  &\texttt{MINRES} \bigl( \boldsymbol{A}_{\beta,\bm{m}}, \boldsymbol{\hat P}^{-1}, \bm{b}, \bm{x}_0 \bigr):
  &
  \norm{\bm{r}_k}_{\boldsymbol{\hat P}^{-1}}
  &= \min_{p\in\mathcal{P}^0_k}
     \norm{p(\boldsymbol{A}_{\beta,\bm{m}}\boldsymbol{\hat P}^{-1})\bm{r}_0}_{\boldsymbol{\hat P}^{-1}},
  \end{aligned}
\end{gather}
\vskip-\belowdisplayskip\vskip\belowdisplayshortskip\noindent
where $\mathcal{P}^0_k$ denotes the set of polynomials of degree at most~$k$ normalized to $p(0)=1$,
$\bm{r}_k=\bm{b}-\boldsymbol{A}_{\beta,\bm{m}}\bm{x}_k$ are the true residuals corresponding to the $k$-th iterates
$\bm{x}_k=[\boldsymbol{\zeta}_k^\top,\,\boldsymbol{\delta}\bm{m}_k^\top]^{\vphantom{\smash[t]{\big|}}\top}$\!\!\!\!,\,
and the norm $\norm{\bm{x}}_{\boldsymbol{M}}=(\bm{x}^\top\boldsymbol{M}\bm{x})^{1/2}$ for a~symmetric positive definite~$\boldsymbol{M}$.
In particular one can see that different residual norms are used.

To assess the complexity of \cref{alg:gnstep-iterative}, we  assume that the black-box preconditioners $\boldsymbol{\hat Q}^{-1}$ and $\boldsymbol{\hat S}^{-1}$ are optimal, i.e., the actions $\boldsymbol{\hat Q}^{-1} \bm{y}_1$, $\boldsymbol{\hat S}^{-1} \bm{y}_2$ on vectors $\bm{y}_1\in{\mathbb R}^K$, $\bm{y}_2\in{\mathbb R}^N$ are performed in $O(K)$ and $O(N)$ floating-point operations, respectively.
In the settings under consideration, we have $M\leq N$ (typically $M\ll N$) and $O(K)=O(N)$.
Moreover, we do not distinguish between complexity for number of floating
point operations and execution times.
A~breakdown of the complexity of the steps in \cref{alg:gnstep-iterative} is as follows:
\begingroup
\setlength\leftmargini{0.75in}
\makeatletter
\interlinepenalty=10000  
\@itempenalty=5000  
\makeatother
\begin{enumerate}
  \item[\cref{ln:H}.]
    $M$ applications of $\boldsymbol{\hat S}^{-1}$, i.e., $O(MN)$;
  \item[\cref{ln:C}.]
    dense matrix-matrix multiply; $O(M^2 N)$;
  \item[\cref{ln:chol}.]
    dense Cholesky factorization; $O(M^3)$;
  \item[\cref{ln:minres}.]
    cost per one MINRES step is $O(MN)$
    because
    $\boldsymbol{Q} \boldsymbol{\zeta}$, $\boldsymbol{D}^\top \boldsymbol{\delta}\bm{m}$, $\boldsymbol{D} \boldsymbol{\zeta}$, $\boldsymbol{\hat Q}^{-1} \bm{y}_1$, and $\boldsymbol{\hat S}^{-1} \bm{y}_2$
    are $O(N)$,
    ${\boldsymbol{J}}_{\bm{m}}^\top ({\boldsymbol{J}}_{\bm{m}}\boldsymbol{\delta}\bm{m})$, $\boldsymbol{\hat H}_{\bm{m}}^\top \bm{y}_2$, and $\boldsymbol{\hat H}_{\bm{m}}{\bm{\mkern 2mu\cdot\mkern 2mu}}$
    are $O(MN)$,
    and $\boldsymbol{L}_{\beta,\bm{m}}^{-\top}{\bm{\mkern 2mu\cdot\mkern 2mu}}$ and $\boldsymbol{L}_{\beta,\bm{m}}^{-1}{\bm{\mkern 2mu\cdot\mkern 2mu}}$
    are $O(M^2)$.
\end{enumerate}
\noindent
\endgroup
If the number of MINRES iterations remains constant independent of $M$ and $N$, one observes that the overall complexity of \cref{alg:gnstep-iterative} is dominated by $O(M^2 N)$ due to the assembly of capacitance matrix on
\cref{ln:C}.
On the other hand, this operation would typically be carried out by the Level~3 BLAS routine \texttt{gemm}, thus very efficiently (in terms of utilizing the theoretical floating point capability of the CPU).
Note that one must not assemble $({\boldsymbol{J}}_{\bm{m}}^\top {\boldsymbol{J}}_{\bm{m}})$, which would be a~dense ${\mathbb R}^{N\times N}$ matrix and thus would
degrade the complexity to $O(N^2)$.
We will demonstrate via the numerical experiments in \cref{sec:experiments}
that the number of MINRES iterations in \cref{alg:gnstep-iterative} tends to be constant.

On the other hand, the simplified \cref{alg:gnstep-iterative-nw} has, by the same reasoning, complexity of only $O(MN)$ per a~MINRES iteration, but the number of MINRES iterations tends to increase as $M$ and $N$
grow, which we will see confirmed in \cref{sec:experiments}.
Moreover, \cref{alg:gnstep-iterative-nw} is not robust w.r.t.~$\beta$; see \cref{fig:timing-beta}.

\section{An application: Electrical resistivity tomography} \label{sec:experiments}

Consider a~conducting medium occupying a~domain $\Omega\subset{\mathbb R}^d$ characterized by an unknown spatially varying electrical conductivity $\sigma_\mathrm{true}\colon\Omega\to(0,\infty)$.
Electrical resistivity tomography (ERT; also known as \emph{the direct current (DC) resistivity method} in the geophysical exploration literature) reconstructs the unknown $\sigma_\mathrm{true}$ from voltage measurements of stationary electric fields excited by known synthetic DC sources.
We model the excitation current by a~source-sink pair of point sources of known DC current strength.
This corresponds physically to a~current source connected to the medium at two distinct points by way of cables (conductors), while the cables themselves are not part of the conductivity model $\sigma_\mathrm{true}$
but are rather represented as a~point source and point sink, respectively.
The response of the medium to this excitation can be measured as a~voltage (potential difference) at two other points in the medium.
By varying the placement of current source/sink and/or the voltage electrode positions one can perform multiple measurements.
Ultimately one wishes to reconstruct a~conductivity distribution $\widetilde\sigma$ which is consistent with these measurements.
A~finite set of such measurements is likely to be explained equally well by multiple different values of~$\widetilde\sigma$, indicating that the problem is underdetermined.
Moreover, the (inverse) problem of reconstructing conductivity from potential measurements is well known to be ill-posed.
As a~selection criterion one can ask for extra smoothness of~$\widetilde\sigma$ and thus regularize the inverse problem.
In any case, it is clear that, except for special cases, it cannot be expected that $\widetilde\sigma=\sigma_\mathrm{true}$.

Consider a~bounded Lipschitz domain $\Omega\subset{\mathbb R}^d$, $d=2,3$ and electrical conductivity $\sigma\in L^\infty(\Omega)$, $\sigma\geq\sigma_0>0$.
Assume ${\partial\Omega}=\overline{\gamma_{\mathrm D}}\cup\overline{\gamma_{\mathrm N}}$ with open and disjoint ${\gamma_{\mathrm D}}$, ${\gamma_{\mathrm N}}$ and such that $\snorm{{\gamma_{\mathrm D}}}>0$.
Note that ${\gamma_{\mathrm D}}$ and ${\gamma_{\mathrm N}}$ are, in general, different from ${\Gamma_{\mathrm D}}$ and ${\Gamma_{\mathrm N}}$ from~\cref{eq:ls}.
We consider the diffusion equation for the stationary electric potential $u$
\begin{subequations} \label{eq_model} \begin{alignat}{2}
    \label{eq_model_om}
    -\operatorname{div} \sigma \nabla u &= f &\qquad &\text{ in } \Omega, \\
    \label{eq_model_gd}
    u &= 0 &\qquad &\text{ on } {\gamma_{\mathrm D}}, \\
    \label{eq_model_gn}
    \tfrac{\partial}{\partial\mathrm{n}}u &= 0 &\qquad &\text{ on } {\gamma_{\mathrm N}},
\end{alignat} \end{subequations}
where we employ the homogeneous boundary conditions~\cref{eq_model_gd,eq_model_gn} for simplicity.
The electric potential ${u_{{x_\mathrm{A}}{x_\mathrm{B}}}}$ for a~unit current source-sink pair ${x_\mathrm{A}} \neq {x_\mathrm{B}}$ in $\Omega\cup{\gamma_{\mathrm N}}$ is then defined as the distributional solution of~\cref{eq_model}
with $f\coloneqq{\delta_{{x_\mathrm{A}}}}-{\delta_{{x_\mathrm{B}}}}$.
Note that it makes sense to place
${x_\mathrm{A}}$ and/or ${x_\mathrm{B}}$ on~${\gamma_{\mathrm N}}$. The distributional solution~${u_{{x_\mathrm{A}}{x_\mathrm{B}}}}$
does not belong to the Sobolev space $H^1(\Omega)$.
Nevertheless, ${u_{{x_\mathrm{A}}{x_\mathrm{B}}}}$ is continuous in $\Omega\setminus({x_\mathrm{A}}\cup{x_\mathrm{B}})$;
see \cite[equation~(3)]{nash-1958}.
One can therefore define the voltage difference
\begin{align} \label{eq_qoi1}
    {u_{{x_\mathrm{A}}{x_\mathrm{B}}}}({x_\mathrm{M}})-{u_{{x_\mathrm{A}}{x_\mathrm{B}}}}({x_\mathrm{N}}) = \langle {\delta_{{x_\mathrm{M}}}}-{\delta_{{x_\mathrm{N}}}}, {u_{{x_\mathrm{A}}{x_\mathrm{B}}}} \rangle
\end{align}
between any two points ${x_\mathrm{M}}$, ${x_\mathrm{N}}\in\Omega\setminus({x_\mathrm{A}}\cup{x_\mathrm{B}})$.
We define the solution operator for equation~\cref{eq_model}:
\begin{align*}
    A^{-1}_\sigma\colon f\mapsto u
    \qquad
    \text{such that } \sigma,\,f,\text{ and }u \text{ satisfy }\cref{eq_model}
    \text{ in the sense of distributions.}
\end{align*}
With this definition we may express the quantity in~\cref{eq_qoi1} as
\begin{equation} \label{eq:ert_qoi}
\begin{aligned}
    \langle {\delta_{{x_\mathrm{M}}}}-{\delta_{{x_\mathrm{N}}}}, A^{-1}_\sigma({\delta_{{x_\mathrm{A}}}}-{\delta_{{x_\mathrm{B}}}}) \rangle
    &=
    \langle {\delta_{{x_\mathrm{A}}}}-{\delta_{{x_\mathrm{B}}}}, A^{-1}_\sigma({\delta_{{x_\mathrm{M}}}}-{\delta_{{x_\mathrm{N}}}}) \rangle \\
    &=
    \int_\Omega \sigma\nabla A^{-1}_\sigma({\delta_{{x_\mathrm{A}}}}-{\delta_{{x_\mathrm{B}}}})\cdot\nabla A^{-1}_\sigma({\delta_{{x_\mathrm{M}}}}-{\delta_{{x_\mathrm{N}}}})
\end{aligned}
\end{equation}
and the G\^{a}teaux derivative of this quantity is readily expressed as\footnote{%
  This follows along the lines of the formula
  $\mathrm{d}(A^{-1}) = -A^{-1}\, \mathrm{d}A\, A^{-1}$, which is valid
  for any invertible matrix~$A$.
  Concerning the G\^{a}teaux derivative of the singular integral~\cref{eq:ert_qoi}
  additional rigor and care in choice of the function spaces is needed,
  but this is out of scope of this work, hence we proceed just formally.
}
\begin{align}
    \label{eq:ert_qoi_derivative}
    \delta\sigma\mapsto
    -\int_\Omega \delta\sigma\nabla A^{-1}_\sigma({\delta_{{x_\mathrm{A}}}}-{\delta_{{x_\mathrm{B}}}})\cdot\nabla A^{-1}_\sigma({\delta_{{x_\mathrm{M}}}}-{\delta_{{x_\mathrm{N}}}}),
\end{align}
which is a~linear functional.

It is convenient to introduce the change of variables $m=\log\sigma$ for the conductivity so that
for $m\in L^\infty(\Omega)$ one has $0 < \exp(\operatorname{ess\,inf}_\Omega m) \allowbreak {}\leq \sigma \leq \exp(\operatorname{ess\,sup}_\Omega m)$.
The solution map $A^{-1}_{\exp(m)}$ is then well defined for all $m\in L^\infty(\Omega)$ as the boundedness condition
$0<\underline\sigma\leq\sigma\leq\overline\sigma<\infty$
is equivalent to $m\in L^\infty(\Omega)$.

A~practical ERT survey consists of multiple measurements using different combinations of points ${x_\mathrm{A}^i}$, ${x_\mathrm{B}^i}$, ${x_\mathrm{M}^i}$,
and~${x_\mathrm{N}^i}$ for $i=1,2,\ldots,M$.
Following~\cref{eq:ert_qoi,eq:ert_qoi_derivative}
we express the quantity of interest and its derivative as:
\begin{multline}
    \label{eq:ert_qoi2}
    \begin{aligned}
      g^i(m) &\coloneqq
      k_i \int_\Omega \exp(m)
      \nabla A^{-1}_{\exp(m)}({\delta_{{x_\mathrm{A}^i}}}-{\delta_{{x_\mathrm{B}^i}}})\cdot\nabla A^{-1}_{\exp(m)}({\delta_{{x_\mathrm{M}^i}}}-{\delta_{{x_\mathrm{N}^i}}}),
      \\
      J^i(m)\delta m &\coloneqq
      -k_i \int_\Omega \delta m \exp(m)
      \nabla A^{-1}_{\exp(m)}({\delta_{{x_\mathrm{A}^i}}}-{\delta_{{x_\mathrm{B}^i}}})\cdot\nabla A^{-1}_{\exp(m)}({\delta_{{x_\mathrm{M}^i}}}-{\delta_{{x_\mathrm{N}^i}}})
    \end{aligned}
    \\
    \text{for } i=1,2,\ldots,M.
\end{multline}
Here we have introduced additional scaling factors $k_i$ given by
\begin{multline}
    \label{eq:ert_qoi_scaling}
    \begin{aligned}
    k_i &\coloneqq \begin{cases}
      \frac{\pi}{-\log\snorm{{x_\mathrm{A}^i}-{x_\mathrm{M}^i}}+\log\snorm{{x_\mathrm{B}^i}-{x_\mathrm{M}^i}}+\log\snorm{{x_\mathrm{A}^i}-{x_\mathrm{N}^i}}-\log\snorm{{x_\mathrm{B}^i}-{x_\mathrm{N}^i}}}
      & d=2, \\[.5em]
      \frac{2\pi}{\snorm{{x_\mathrm{A}^i}-{x_\mathrm{M}^i}}^{-1}-\snorm{{x_\mathrm{B}^i}-{x_\mathrm{M}^i}}^{-1}-\snorm{{x_\mathrm{A}^i}-{x_\mathrm{N}^i}}^{-1}+\snorm{{x_\mathrm{B}^i}-{x_\mathrm{N}^i}}^{-1}}
      & d=3,
    \end{cases}
    \end{aligned}
    \\
    \text{for } i=1,2,\ldots,M.
\end{multline}
These \emph{geometric factors} only depend on the coordinates of the electrodes.
Their purpose is that the original voltage measurement~\cref{eq:ert_qoi}
is transformed into a~quantity known as \emph{apparent resistivity}\footnote{%
  A~measurement $g^i(m)$ from~\cref{eq:ert_qoi2} gives apparent constant
  resistivity of homogeneous half space. Precisely, it holds true
  that $g^i(\log\sigma_0) = 1/\sigma_0$, for a~constant $\sigma_0>0$,
  half-space domain $\Omega=\{x\in{\mathbb R}^d,\,x_d>0\}$,
  and ${x_\mathrm{A}^i},\,{x_\mathrm{B}^i},\,{x_\mathrm{M}^i},\,{x_\mathrm{N}^i} \in \{x\in{\mathbb R}^d,\,x_d=0\}$.
  This is derived using Green's functions for the Laplace Dirichlet
  problem in half space.
}.
This is a~commonly applied method of scaling the data $g^i(m)-{g^i_{\mathrm{obs}}}$, $i=1,\ldots,M$.

\begin{figure}[t]
  \centering
  \includegraphics[width=\textwidth,trim=69 40 77 9,clip]{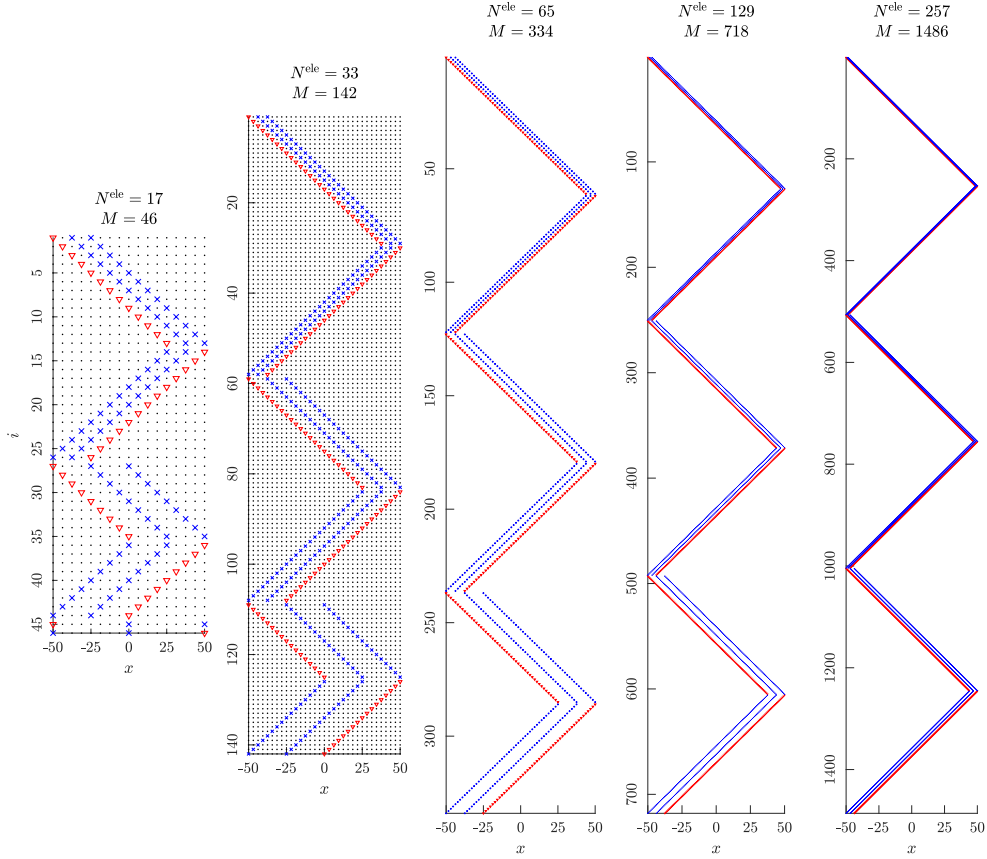}
  \caption{%
    The first five \emph{pole-dipole}
    electrode configurations used in the numerical examples,
    starting from $N^\mathrm{ele}=17$ distinct electrode positions (left)
    and reaching as many as $N^\mathrm{ele}=257$ positions (right).
    The vertical axes enumerate measurement number $i=1,2,\ldots,M$
    and the horizontal axes represent electrode $x$-position
    with $N^\mathrm{ele}$ distinct equidistant positions in the
    interval $x\in[-50,50]$.
    Transmitter electrodes ${x_\mathrm{A}^i}$ ($\color{red}\triangledown$),
    receiver electrodes ${x_\mathrm{M}^i}$ ($\color{blue}\times$),
    ${x_\mathrm{N}^i}$ ($\color{blue}\times$),
    second transmitter electrode ${x_\mathrm{B}^i}$ placed at~$\infty$
    (not shown),
    $i=1,2,\ldots,M$.
    The measurements $i=1,2,\ldots,2N^\mathrm{ele}-8$
    use the spacing $2$, e.g., for $i=1$, the electrodes are
    at positions $1,3,5$, for $i=2$ the positions $2,4,6$, etc.
    The measurements $i=2N^\mathrm{ele}-7,\ldots,
    4N^\mathrm{ele}-24$ use the spacing $4$
    and
    the measurements $i=4N^\mathrm{ele}-23,\ldots,
    6N^\mathrm{ele}-56$ use the spacing $8$.
    Total number of measurements is $M=6N^\mathrm{ele}-56$.
  }
  \label{fig:survey}
\end{figure}
In practice, the placement of electrodes ${x_\mathrm{A}^i}$, ${x_\mathrm{B}^i}$,
${x_\mathrm{M}^i}$, and ${x_\mathrm{N}^i}$ is critical for the goal of approximating
the original conductivity distribution, i.e., $\widetilde\sigma
\approx\sigma_\mathrm{true}$.
The geophysics literature contains a~number of established electrode
placement designs; see, e.g., \cite[section~8.5]{TelfordEtAl1990}
or~\cite{uhlemann2018} and the references therein.
In the examples below we consider what is known as a~\emph{pole-dipole configuration}.
\Cref{fig:survey} shows a~sequence of one-dimensional electrode configurations, which exhibit increasing measurement resolution, but with sensitivity only in regions increasingly closer to the surface as the
configuration is refined.
This sequence was chosen to obtain a~series of problems which are each meaningful for the underlying inverse problem and at the same time illustrate the performance of the preconditioners across a~wide range of values for the finite element mesh size, number of measurements, and regularization parameter.
In the framework of~\cref{eq:ert_qoi2,eq:ert_qoi_scaling}, the ${x_\mathrm{B}^i}$-electrode is modeled as an electrode placed at~$\infty$, resp. at~${\gamma_{\mathrm D}}$ in the context of the boundary datum~\cref{eq_model_gd}.
Hence ${\delta_{{x_\mathrm{B}^i}}}$ does not contribute to~\cref{eq:ert_qoi2} and the
factors~$k_i$ are obtained by taking limit \smash{$\snorm{{x_\mathrm{B}^i}}\to\infty$}
in~\cref{eq:ert_qoi_scaling}.
These one-dimensional configurations are typically used in ERT surveys along the upper boundary of a~two-dimensional vertical cross section.
For surveys over a~three-dimensional region, it is common to construct a~two-dimensional surface electrode configuration as the Cartesian product of the one-dimensional pattern.

\begin{figure}[p]
  \centering
  \includegraphics[width=\textwidth]{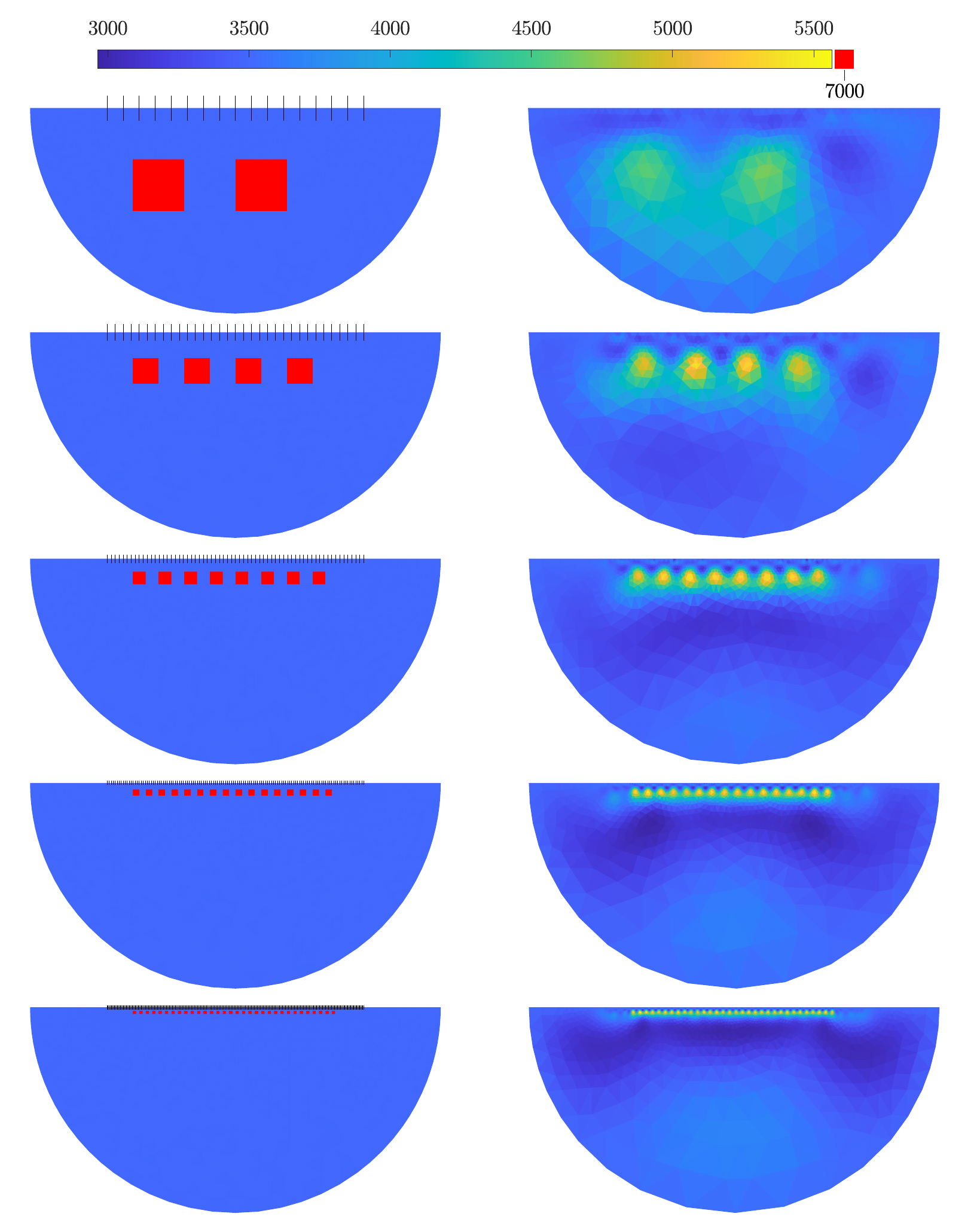}
  \unskip
  \caption{%
    True resistivities $\frac{1}{\sigma_\mathrm{true}}$ (on the left)
    and the result of inversion using \cref{alg:gnstep-iterative}
    (on the right) for the 2D example.
    The series of configurations (from top to bottom) corresponds
    to the first five configurations of electrodes; see \cref{fig:survey}.
    Majority of the medium has background resistivity
    $\frac{1}{\sigma_\mathrm{ref}}=3500$ (blue, on the left)
    with presence of anomaly resistivity $7000$ (red, on the left).
    The positions of the electrodes at the surface
    are indicated by the black vertical bar ($|$, on the left).
  }
  \label{fig:ex2d}
\end{figure}

In the following we illustrate the performance of the aforementioned algorithms with a~sequence of parameter identification experiments in an idealized ERT setting.
We aim to reconstruct a~priori known conductivity anomaly against a~homogeneous background.
We consider a~sequence of problems involving a~checkerboard anomaly structure of increasing complexity with decreasing depth in accordance with the sensitivity and resolution capability of the chosen electrode configuration designs.
Sequences of problems in two (the left column in \cref{fig:ex2d}) and three spatial dimensions
(the top row in \cref{fig:ex3d}) are carefully chosen to work well with the aforementioned electrode configuration. In particular, because the spacing of electrodes decreases with finer configurations,
the survey is only sensitive in an increasingly shallow region below the surface.
This rather artificial scenario allows us to reconstruct an~increasingly finer pattern
with only $M=O(N)$ measurements (see \cref{fig:survey}), thus allowing us to increase
the parameters~$M$ and~$N$ many times before exhausting compute resources
(see \cref{tab:experiments}).
\begin{table}[t]
  \capstart
  \newcolumntype{H}{>{\setbox0=\hbox\bgroup}c<{\egroup}@{}}  
  \newcolumntype{Z}{>{\setbox0=\hbox\bgroup}c<{\egroup}@{\hspace*{-\tabcolsep}}}  
  \def\dncsign{$\smash{\dagger}$}
  \def\dnc{\rlap{\dncsign}}
  \caption{%
    Performance characteristics of the numerical experiments.
    Timings $t_{\labelcref{ln:H}}$, $t_{\labelcref{ln:C}}$, and
    $t_{\labelcref{ln:chol}}$ for substeps of
    \cref{alg:gnstep-iterative} for each Gauss--Newton step~$i$.
    Number of MINRES iterations $n_\mathrm{iter}$ and overall
    runtime $t_\mathrm{norm}$ for
    \cref{alg:gnstep-iterative,alg:gnstep-iterative-nw}
    to solve the normal equations $\boldsymbol{A}_{\beta,\bm{m}}\bm{x}=\bm{b}$
    within tolerance $\norm{\bm{r}_k}_2/\norm{\bm{b}}_2 \leq 10^{-7}$
    in the Euclidean norm.
    Cases marked~\dncsign{} did not converge to the
    prescribed tolerance in $2(K+N)$ iterations
    (recall that $\boldsymbol{A}_{\beta,\bm{m}}\in{\mathbb R}^{(K+N)\times(K+N)}$);
    tolerance $10^{-6}$ was reached in all these cases.
  }
  \centering
  \setlength{\tabcolsep}{5.6pt}
  \scriptsize%
  \renewcommand\RSsmallest{4.5pt}
  \begin{tabular}{cZ}
    & \vphantom{Alg.1} \\
    & \vphantom{$N^\mathrm{ele}$} \\ \hline
    \multirow{18}{*}{\begin{adjustbox}{angle=90}2D ($\beta=0.1$)\end{adjustbox}}
    & 1 \\ & 1 \\ & 1 \\ & 1 \\ & 1 \\ & 1 \\ & 1 \\ & 1 \\ & 1 \\
    & 1 \\ & 1 \\ & 1 \\ & 1 \\ & 1 \\ & 1 \\ & 1 \\ & 1 \\ & 1 \\
    \hline
    \multirow{10}{*}{\begin{adjustbox}{angle=90}3D ($\beta=10^5$)\end{adjustbox}}
    & 1 \\ & 1 \\ & 1 \\ & 1 \\ & 1 \\ & 1 \\
    & 1 \\ & 1 \\ & 1 \\ & 1 \\
  \end{tabular}
  \hspace*{-\floatsep}%
  \begin{tabular}{HHrrrr|rrrrr}
    \multirow{2}{*}{$n_x$}
    & \multirow{2}{*}{$n_y$}
    & \multirow{2}{*}{$N^\mathrm{ele}$}
    & \multirow{2}{*}{$N$}
    & \multirow{2}{*}{$M$}
    & \multirow{2}{*}{$i$}
    & \multicolumn{5}{c}{\cref{alg:gnstep-iterative}}
    \\
    & & & & &
    & \multicolumn{1}{c}{$t_{\labelcref{ln:H}}\;[\mathrm{s}]$}
    & \multicolumn{1}{c}{$t_{\labelcref{ln:C}}\;[\mathrm{s}]$}
    & \multicolumn{1}{c}{$t_{\labelcref{ln:chol}}\;[\mathrm{s}]$}
    & \multicolumn{1}{c}{$n_\mathrm{iter}$}
    & \multicolumn{1}{c}{$t_\mathrm{norm}\;[\mathrm{s}]$}
    \\ \hline
\multirow[t]{2}{*}{4} & \multirow[t]{18}{*}{1} & \multirow[t]{2}{*}{17} & \multirow[t]{2}{*}{840} & \multirow[t]{2}{*}{46} & 1 & 0.01 & 0.0019 & 0.0001 & 4 & 0.02\\ &  &  &  &  & 2 & 0.01 & 0.0007 & 0.0001 & 11 & 0.02\\\multirow[t]{2}{*}{8} &  & \multirow[t]{2}{*}{33} & \multirow[t]{2}{*}{1584} & \multirow[t]{2}{*}{142} & 1 & 0.03 & 0.0029 & 0.0050 & 4 & 0.07\\ &  &  &  &  & 2 & 0.03 & 0.0019 & 0.0024 & 12 & 0.05\\\multirow[t]{2}{*}{16} &  & \multirow[t]{2}{*}{65} & \multirow[t]{2}{*}{3140} & \multirow[t]{2}{*}{334} & 1 & 0.16 & 0.0116 & 0.0010 & 4 & 0.19\\ &  &  &  &  & 2 & 0.15 & 0.0113 & 0.0007 & 14 & 0.22\\\multirow[t]{2}{*}{32} &  & \multirow[t]{2}{*}{129} & \multirow[t]{2}{*}{6012} & \multirow[t]{2}{*}{718} & 1 & 0.61 & 0.0930 & 0.0026 & 4 & 0.80\\ &  &  &  &  & 2 & 0.63 & 0.0885 & 0.0118 & 14 & 1.09\\\multirow[t]{2}{*}{64} &  & \multirow[t]{2}{*}{257} & \multirow[t]{2}{*}{11644} & \multirow[t]{2}{*}{1486} & 1 & 2.49 & 0.4523 & 0.0133 & 4 & 3.24\\ &  &  &  &  & 2 & 2.64 & 0.4656 & 0.0139 & 17 & 3.88\\\multirow[t]{2}{*}{128} &  & \multirow[t]{2}{*}{513} & \multirow[t]{2}{*}{22884} & \multirow[t]{2}{*}{3022} & 1 & 11.01 & 1.8966 & 0.0668 & 4 & 13.98\\ &  &  &  &  & 2 & 10.78 & 1.9269 & 0.0567 & 14 & 15.40\\\multirow[t]{2}{*}{256} &  & \multirow[t]{2}{*}{1025} & \multirow[t]{2}{*}{44848} & \multirow[t]{2}{*}{6094} & 1 & 44.22 & 14.6608 & 0.3928 & 4 & 63.39\\ &  &  &  &  & 2 & 44.93 & 14.6799 & 0.3932 & 12 & 69.61\\\multirow[t]{2}{*}{512} &  & \multirow[t]{2}{*}{2049} & \multirow[t]{2}{*}{89608} & \multirow[t]{2}{*}{12238} & 1 & 190.69 & 117.3031 & 3.0011 & 6 & 330.03\\ &  &  &  &  & 2 & 199.42 & 117.3992 & 2.9831 & 14 & 357.53\\\multirow[t]{2}{*}{1024} &  & \multirow[t]{2}{*}{4097} & \multirow[t]{2}{*}{178232} & \multirow[t]{2}{*}{24526} & 1 & 848.10 & 932.6321 & 23.4441 & 7 & 1943.20\\ &  &  &  &  & 2 & 847.55 & 932.3276 & 22.3924 & 14 & 1963.20 %
\unskip
    \\ \hline
\multirow[t]{10}{*}{2} & \multirow[t]{2}{*}{2} & \multirow[t]{2}{*}{81} & \multirow[t]{2}{*}{120192} & \multirow[t]{2}{*}{216} & 1 & 5.58 & 0.1238 & 0.0005 & 130 & 24.19\\ &  &  &  &  & 2 & 5.61 & 0.1209 & 0.0004 & 132 & 20.98\\ & \multirow[t]{2}{*}{3} & \multirow[t]{2}{*}{169} & \multirow[t]{2}{*}{262464} & \multirow[t]{2}{*}{728} & 1 & 54.31 & 2.6834 & 0.0036 & 123 & 139.58\\ &  &  &  &  & 2 & 53.76 & 2.6758 & 0.0028 & 125 & 169.42\\ & \multirow[t]{2}{*}{4} & \multirow[t]{2}{*}{289} & \multirow[t]{2}{*}{452736} & \multirow[t]{2}{*}{1564} & 1 & 171.92 & 10.1975 & 0.0155 & 124 & 394.06\\ &  &  &  &  & 2 & 176.60 & 10.2716 & 0.0153 & 126 & 392.27\\ & \multirow[t]{2}{*}{5} & \multirow[t]{2}{*}{441} & \multirow[t]{2}{*}{679296} & \multirow[t]{2}{*}{2940} & 1 & 536.86 & 53.4599 & 0.1874 & 136 & 1208.20\\ &  &  &  &  & 2 & 516.93 & 53.4413 & 0.0981 & 139 & 1259.60\\ & \multirow[t]{2}{*}{6} & \multirow[t]{2}{*}{625} & \multirow[t]{2}{*}{937408} & \multirow[t]{2}{*}{4700} & 1 & 1146.50 & 184.4371 & 0.2116 & 131 & 2840.50\\ &  &  &  &  & 2 & 1270.10 & 183.5881 & 0.2154 & 136 & 2988.60 %
\unskip
  \end{tabular}
  \hspace*{-\floatsep}%
  \begin{tabular}{HHHHHHHHH|rr}
    &&&&&&&&
    & \multicolumn{2}{c}{\cref{alg:gnstep-iterative-nw}}
    \\
    &&&&&&&&
    & \multicolumn{1}{c}{$n_\mathrm{iter}$}
    & \multicolumn{1}{c}{$t_\mathrm{norm}\;[\mathrm{s}]$}
    \\ \hline
\multirow[t]{2}{*}{4} & \multirow[t]{10}{*}{1} & \multirow[t]{2}{*}{17} & \multirow[t]{2}{*}{840} & \multirow[t]{2}{*}{46} & 1 & 0.00 & 0.0000 & 0.0000 & 79 & 0.02\\ &  &  &  &  & 2 & 0.00 & 0.0000 & 0.0000 & 286 & 0.07\\\multirow[t]{2}{*}{8} &  & \multirow[t]{2}{*}{33} & \multirow[t]{2}{*}{1584} & \multirow[t]{2}{*}{142} & 1 & 0.00 & 0.0000 & 0.0000 & 737 & 0.38\\ &  &  &  &  & 2 & 0.00 & 0.0000 & 0.0000 & 3526 & 1.58\\\multirow[t]{2}{*}{16} &  & \multirow[t]{2}{*}{65} & \multirow[t]{2}{*}{3140} & \multirow[t]{2}{*}{334} & 1 & 0.00 & 0.0000 & 0.0000 & 1921 & 6.46\\ &  &  &  &  & 2 & 0.00 & 0.0000 & 0.0000 & 15878\dnc & 39.69\\\multirow[t]{2}{*}{32} &  & \multirow[t]{2}{*}{129} & \multirow[t]{2}{*}{6012} & \multirow[t]{2}{*}{718} & 1 & 0.00 & 0.0000 & 0.0000 & 3441 & 33.71\\ &  &  &  &  & 2 & 0.00 & 0.0000 & 0.0000 & 29883\dnc & 331.63\\\multirow[t]{2}{*}{64} &  & \multirow[t]{2}{*}{257} & \multirow[t]{2}{*}{11644} & \multirow[t]{2}{*}{1486} & 1 & 0.00 & 0.0000 & 0.0000 & 5700 & 180.17\\ &  &  &  &  & 2 & 0.00 & 0.0000 & 0.0000 & 58798\dnc & 2016.60
\unskip
    \\ &&&&&&&&&&
    \\ &&&&&&&&&&
    \\ &&&&&&&&&&
    \\ &&&&&&&&&&
    \\ &&&&&&&&&&
    \\ &&&&&&&&&&
    \\ &&&&&&&&&&
    \\ &&&&&&&&&&
    \\ \hline
\multirow[t]{10}{*}{2} & \multirow[t]{2}{*}{2} & \multirow[t]{2}{*}{81} & \multirow[t]{2}{*}{120192} & \multirow[t]{2}{*}{216} & 1 & 0.00 & 0.0000 & 0.0000 & 156 & 15.77\\ &  &  &  &  & 2 & 0.00 & 0.0000 & 0.0000 & 160 & 15.43\\ & \multirow[t]{2}{*}{3} & \multirow[t]{2}{*}{169} & \multirow[t]{2}{*}{262464} & \multirow[t]{2}{*}{728} & 1 & 0.00 & 0.0000 & 0.0000 & 210 & 94.94\\ &  &  &  &  & 2 & 0.00 & 0.0000 & 0.0000 & 215 & 72.15\\ & \multirow[t]{2}{*}{4} & \multirow[t]{2}{*}{289} & \multirow[t]{2}{*}{452736} & \multirow[t]{2}{*}{1564} & 1 & 0.00 & 0.0000 & 0.0000 & 291 & 277.57\\ &  &  &  &  & 2 & 0.00 & 0.0000 & 0.0000 & 297 & 314.99\\ & \multirow[t]{2}{*}{5} & \multirow[t]{2}{*}{441} & \multirow[t]{2}{*}{679296} & \multirow[t]{2}{*}{2940} & 1 & 0.00 & 0.0000 & 0.0000 & 425 & 1363.20\\ &  &  &  &  & 2 & 0.00 & 0.0000 & 0.0000 & 438 & 1204.00\\ & \multirow[t]{2}{*}{6} & \multirow[t]{2}{*}{625} & \multirow[t]{2}{*}{937408} & \multirow[t]{2}{*}{4700} & 1 & 0.00 & 0.0000 & 0.0000 & 526 & 2989.20\\ &  &  &  &  & 2 & 0.00 & 0.0000 & 0.0000 & 542 & 3100.70 %
\unskip
  \end{tabular}
  \label{tab:experiments}
\end{table}

\begin{figure}[p]
  \centering
  \footnotesize
  \def\svgwidth{1.\textwidth}
  %
\begingroup%
  \makeatletter%
  \providecommand\color[2][]{%
    \errmessage{(Inkscape) Color is used for the text in Inkscape, but the package 'color.sty' is not loaded}%
    \renewcommand\color[2][]{}%
  }%
  \providecommand\transparent[1]{%
    \errmessage{(Inkscape) Transparency is used (non-zero) for the text in Inkscape, but the package 'transparent.sty' is not loaded}%
    \renewcommand\transparent[1]{}%
  }%
  \providecommand\rotatebox[2]{#2}%
  \newcommand*\fsize{\dimexpr\f@size pt\relax}%
  \newcommand*\lineheight[1]{\fontsize{\fsize}{#1\fsize}\selectfont}%
  \ifx\svgwidth\undefined%
    \setlength{\unitlength}{1538.00003052bp}%
    \ifx\svgscale\undefined%
      \relax%
    \else%
      \setlength{\unitlength}{\unitlength * \real{\svgscale}}%
    \fi%
  \else%
    \setlength{\unitlength}{\svgwidth}%
  \fi%
  \global\let\svgwidth\undefined%
  \global\let\svgscale\undefined%
  \makeatother%
  \begin{picture}(1,1.25422627)%
    \lineheight{1}%
    \setlength\tabcolsep{0pt}%
    \put(0,0){\includegraphics[width=\unitlength,page=1]{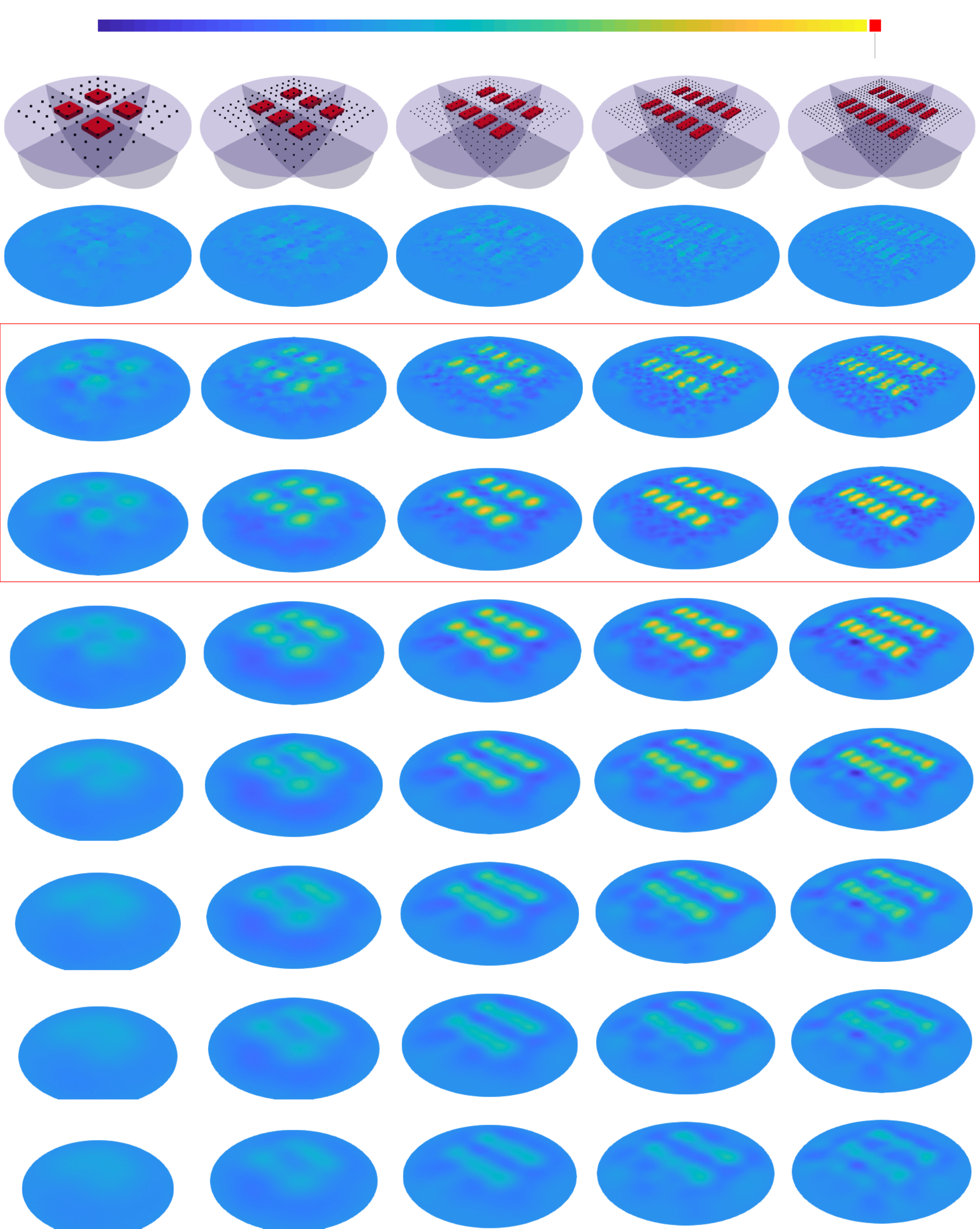}}%
    \put(0.1,1.23827932){\color[rgb]{0,0,0}\makebox(0,0)[t]{\lineheight{1.25}\smash{\begin{tabular}[t]{c}$3203$\end{tabular}}}}%
    \put(0.2960338,1.23827932){\color[rgb]{0,0,0}\makebox(0,0)[t]{\lineheight{1.25}\smash{\begin{tabular}[t]{c}$3435$\end{tabular}}}}%
    \put(0.49206761,1.23827932){\color[rgb]{0,0,0}\makebox(0,0)[t]{\lineheight{1.25}\smash{\begin{tabular}[t]{c}$3667$\end{tabular}}}}%
    \put(0.68810146,1.23827932){\color[rgb]{0,0,0}\makebox(0,0)[t]{\lineheight{1.25}\smash{\begin{tabular}[t]{c}$3898$\end{tabular}}}}%
    \put(0.88413527,1.23827932){\color[rgb]{0,0,0}\makebox(0,0)[t]{\lineheight{1.25}\smash{\begin{tabular}[t]{c}$4130$\end{tabular}}}}%
    \put(0.893238,1.18106216){\color[rgb]{0,0,0}\makebox(0,0)[t]{\lineheight{1.25}\smash{\begin{tabular}[t]{c}$7000$\end{tabular}}}}%
    \put(0.00520156,0.9313186){\color[rgb]{0,0,0}\makebox(0,0)[lt]{\lineheight{1.25}\smash{\begin{tabular}[t]{l}$z = 0.00$\end{tabular}}}}%
    \put(0.20520156,0.9313186){\color[rgb]{0,0,0}\makebox(0,0)[lt]{\lineheight{1.25}\smash{\begin{tabular}[t]{l}$z = 0.00$\end{tabular}}}}%
    \put(0.40520157,0.9313186){\color[rgb]{0,0,0}\makebox(0,0)[lt]{\lineheight{1.25}\smash{\begin{tabular}[t]{l}$z = 0.00$\end{tabular}}}}%
    \put(0.60520155,0.9313186){\color[rgb]{0,0,0}\makebox(0,0)[lt]{\lineheight{1.25}\smash{\begin{tabular}[t]{l}$z = 0.00$\end{tabular}}}}%
    \put(0.80520158,0.9313186){\color[rgb]{0,0,0}\makebox(0,0)[lt]{\lineheight{1.25}\smash{\begin{tabular}[t]{l}$z = 0.00$\end{tabular}}}}%
    \put(0.00520156,0.79929517){\color[rgb]{0,0,0}\makebox(0,0)[lt]{\lineheight{1.25}\smash{\begin{tabular}[t]{l}$z = 5.00$\end{tabular}}}}%
    \put(0.20520156,0.79929517){\color[rgb]{0,0,0}\makebox(0,0)[lt]{\lineheight{1.25}\smash{\begin{tabular}[t]{l}$z = 3.33$\end{tabular}}}}%
    \put(0.40520157,0.79929517){\color[rgb]{0,0,0}\makebox(0,0)[lt]{\lineheight{1.25}\smash{\begin{tabular}[t]{l}$z = 2.50$\end{tabular}}}}%
    \put(0.60520155,0.79929517){\color[rgb]{0,0,0}\makebox(0,0)[lt]{\lineheight{1.25}\smash{\begin{tabular}[t]{l}$z = 2.00$\end{tabular}}}}%
    \put(0.80520158,0.79929517){\color[rgb]{0,0,0}\makebox(0,0)[lt]{\lineheight{1.25}\smash{\begin{tabular}[t]{l}$z = 1.67$\end{tabular}}}}%
    \put(0.00520156,0.66727113){\color[rgb]{0,0,0}\makebox(0,0)[lt]{\lineheight{1.25}\smash{\begin{tabular}[t]{l}$z = 10.00$\end{tabular}}}}%
    \put(0.20520156,0.66727113){\color[rgb]{0,0,0}\makebox(0,0)[lt]{\lineheight{1.25}\smash{\begin{tabular}[t]{l}$z = 6.67$\end{tabular}}}}%
    \put(0.40520157,0.66727113){\color[rgb]{0,0,0}\makebox(0,0)[lt]{\lineheight{1.25}\smash{\begin{tabular}[t]{l}$z = 5.00$\end{tabular}}}}%
    \put(0.60520155,0.66727113){\color[rgb]{0,0,0}\makebox(0,0)[lt]{\lineheight{1.25}\smash{\begin{tabular}[t]{l}$z = 4.00$\end{tabular}}}}%
    \put(0.80520158,0.66727113){\color[rgb]{0,0,0}\makebox(0,0)[lt]{\lineheight{1.25}\smash{\begin{tabular}[t]{l}$z = 3.33$\end{tabular}}}}%
    \put(0.00520156,0.53524705){\color[rgb]{0,0,0}\makebox(0,0)[lt]{\lineheight{1.25}\smash{\begin{tabular}[t]{l}$z = 15.00$\end{tabular}}}}%
    \put(0.20520156,0.53524705){\color[rgb]{0,0,0}\makebox(0,0)[lt]{\lineheight{1.25}\smash{\begin{tabular}[t]{l}$z = 10.00$\end{tabular}}}}%
    \put(0.40520157,0.53524705){\color[rgb]{0,0,0}\makebox(0,0)[lt]{\lineheight{1.25}\smash{\begin{tabular}[t]{l}$z = 7.50$\end{tabular}}}}%
    \put(0.60520155,0.53524705){\color[rgb]{0,0,0}\makebox(0,0)[lt]{\lineheight{1.25}\smash{\begin{tabular}[t]{l}$z = 6.00$\end{tabular}}}}%
    \put(0.80520158,0.53524705){\color[rgb]{0,0,0}\makebox(0,0)[lt]{\lineheight{1.25}\smash{\begin{tabular}[t]{l}$z = 5.00$\end{tabular}}}}%
    \put(0.00520156,0.40322499){\color[rgb]{0,0,0}\makebox(0,0)[lt]{\lineheight{1.25}\smash{\begin{tabular}[t]{l}$z = 20.00$\end{tabular}}}}%
    \put(0.20520156,0.40322499){\color[rgb]{0,0,0}\makebox(0,0)[lt]{\lineheight{1.25}\smash{\begin{tabular}[t]{l}$z = 13.33$\end{tabular}}}}%
    \put(0.40520157,0.40322499){\color[rgb]{0,0,0}\makebox(0,0)[lt]{\lineheight{1.25}\smash{\begin{tabular}[t]{l}$z = 10.00$\end{tabular}}}}%
    \put(0.60520155,0.40322499){\color[rgb]{0,0,0}\makebox(0,0)[lt]{\lineheight{1.25}\smash{\begin{tabular}[t]{l}$z = 8.00$\end{tabular}}}}%
    \put(0.80520158,0.40322499){\color[rgb]{0,0,0}\makebox(0,0)[lt]{\lineheight{1.25}\smash{\begin{tabular}[t]{l}$z = 6.67$\end{tabular}}}}%
    \put(0.00520156,0.27120285){\color[rgb]{0,0,0}\makebox(0,0)[lt]{\lineheight{1.25}\smash{\begin{tabular}[t]{l}$z = 25.00$\end{tabular}}}}%
    \put(0.20520156,0.27120285){\color[rgb]{0,0,0}\makebox(0,0)[lt]{\lineheight{1.25}\smash{\begin{tabular}[t]{l}$z = 16.67$\end{tabular}}}}%
    \put(0.40520157,0.27120285){\color[rgb]{0,0,0}\makebox(0,0)[lt]{\lineheight{1.25}\smash{\begin{tabular}[t]{l}$z = 12.50$\end{tabular}}}}%
    \put(0.60520155,0.27120285){\color[rgb]{0,0,0}\makebox(0,0)[lt]{\lineheight{1.25}\smash{\begin{tabular}[t]{l}$z = 10.00$\end{tabular}}}}%
    \put(0.80520158,0.27120285){\color[rgb]{0,0,0}\makebox(0,0)[lt]{\lineheight{1.25}\smash{\begin{tabular}[t]{l}$z = 8.33$\end{tabular}}}}%
    \put(0.00520156,0.13917428){\color[rgb]{0,0,0}\makebox(0,0)[lt]{\lineheight{1.25}\smash{\begin{tabular}[t]{l}$z = 30.00$\end{tabular}}}}%
    \put(0.20520156,0.13917428){\color[rgb]{0,0,0}\makebox(0,0)[lt]{\lineheight{1.25}\smash{\begin{tabular}[t]{l}$z = 20.00$\end{tabular}}}}%
    \put(0.40520157,0.13917428){\color[rgb]{0,0,0}\makebox(0,0)[lt]{\lineheight{1.25}\smash{\begin{tabular}[t]{l}$z = 15.00$\end{tabular}}}}%
    \put(0.60520155,0.13917428){\color[rgb]{0,0,0}\makebox(0,0)[lt]{\lineheight{1.25}\smash{\begin{tabular}[t]{l}$z = 12.00$\end{tabular}}}}%
    \put(0.80520158,0.13917428){\color[rgb]{0,0,0}\makebox(0,0)[lt]{\lineheight{1.25}\smash{\begin{tabular}[t]{l}$z = 10.00$\end{tabular}}}}%
    \put(0.00520156,0.00715215){\color[rgb]{0,0,0}\makebox(0,0)[lt]{\lineheight{1.25}\smash{\begin{tabular}[t]{l}$z = 35.00$\end{tabular}}}}%
    \put(0.20520156,0.00715215){\color[rgb]{0,0,0}\makebox(0,0)[lt]{\lineheight{1.25}\smash{\begin{tabular}[t]{l}$z = 23.33$\end{tabular}}}}%
    \put(0.40520157,0.00715215){\color[rgb]{0,0,0}\makebox(0,0)[lt]{\lineheight{1.25}\smash{\begin{tabular}[t]{l}$z = 17.50$\end{tabular}}}}%
    \put(0.60520155,0.00715215){\color[rgb]{0,0,0}\makebox(0,0)[lt]{\lineheight{1.25}\smash{\begin{tabular}[t]{l}$z = 14.00$\end{tabular}}}}%
    \put(0.80520158,0.00715215){\color[rgb]{0,0,0}\makebox(0,0)[lt]{\lineheight{1.25}\smash{\begin{tabular}[t]{l}$z = 11.67$\end{tabular}}}}%
  \end{picture}%
\endgroup%
   \unskip
  \caption{%
    Series of 3D computational examples of increasing difficulty
    (left to right). Domain indicated
    by slices $x=0$, $y=0$, and $z=0$ (top row); electrode positions
    indicated by black dots (top row).
    True resistivities $\frac{1}{\sigma_\mathrm{true}}$ in the
    majority of the medium is $\frac{1}{\sigma_\mathrm{ref}}=3500$,
    with presence of anomaly resistivity $7000$ (top row, red blocks).
    Result of inversion using \cref{alg:gnstep-iterative}
    with $\beta=10^5$
    (remaining rows; sections through $z=\mathrm{const}$ planes);
    the slices framed in the red frame correspond to the top and
    the bottom of the anomalous resistivity (red blocks).
  }
  \label{fig:ex3d}
\end{figure}
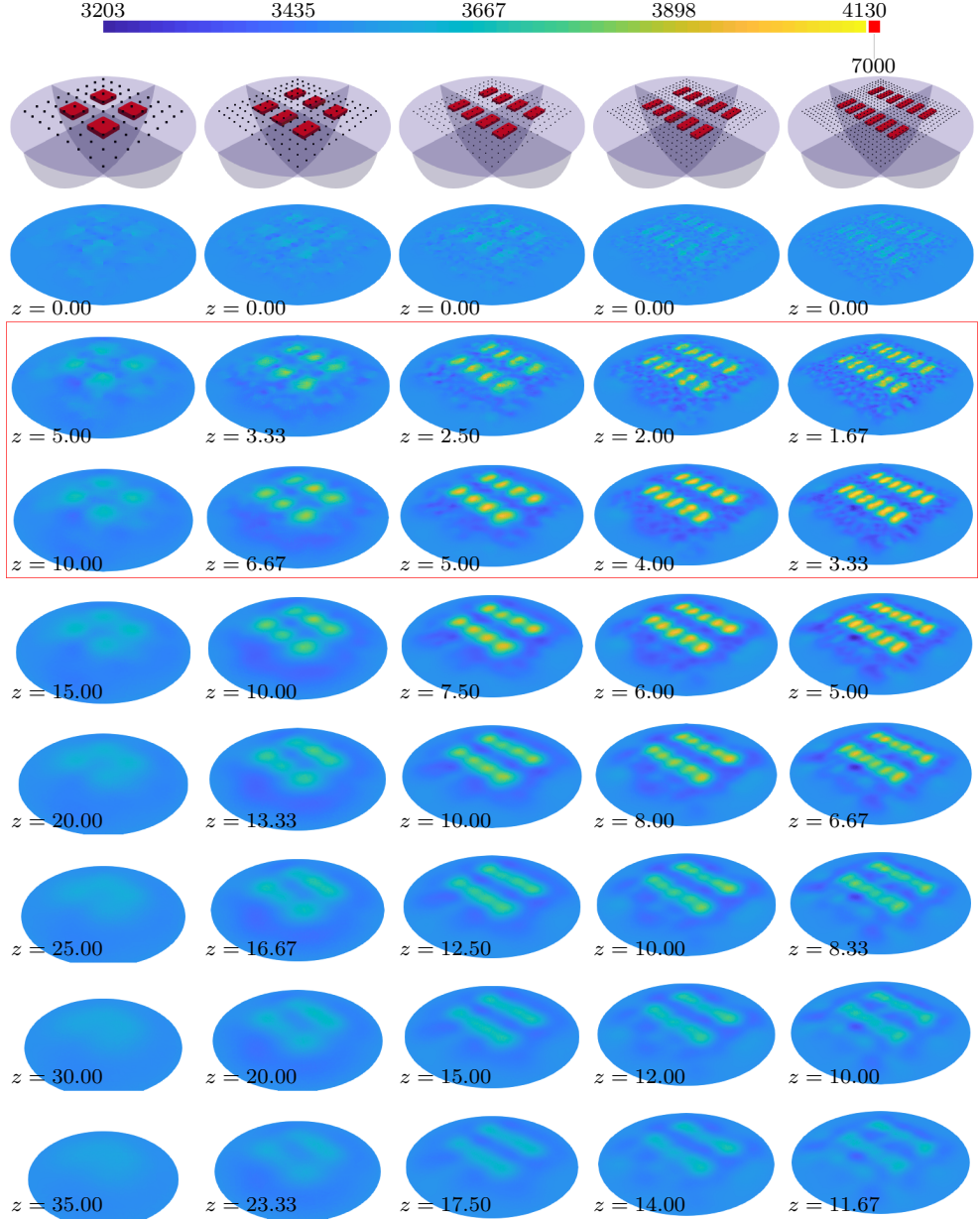

\begin{figure}[p]
  \centering
  \footnotesize
  \def\svgwidth{1.\textwidth}
  %
\begingroup%
  \makeatletter%
  \providecommand\color[2][]{%
    \errmessage{(Inkscape) Color is used for the text in Inkscape, but the package 'color.sty' is not loaded}%
    \renewcommand\color[2][]{}%
  }%
  \providecommand\transparent[1]{%
    \errmessage{(Inkscape) Transparency is used (non-zero) for the text in Inkscape, but the package 'transparent.sty' is not loaded}%
    \renewcommand\transparent[1]{}%
  }%
  \providecommand\rotatebox[2]{#2}%
  \newcommand*\fsize{\dimexpr\f@size pt\relax}%
  \newcommand*\lineheight[1]{\fontsize{\fsize}{#1\fsize}\selectfont}%
  \ifx\svgwidth\undefined%
    \setlength{\unitlength}{1538.00005399bp}%
    \ifx\svgscale\undefined%
      \relax%
    \else%
      \setlength{\unitlength}{\unitlength * \real{\svgscale}}%
    \fi%
  \else%
    \setlength{\unitlength}{\svgwidth}%
  \fi%
  \global\let\svgwidth\undefined%
  \global\let\svgscale\undefined%
  \makeatother%
  \begin{picture}(1,1.30078615)%
    \lineheight{1}%
    \setlength\tabcolsep{0pt}%
    \put(0,0){\includegraphics[width=\unitlength,page=1]{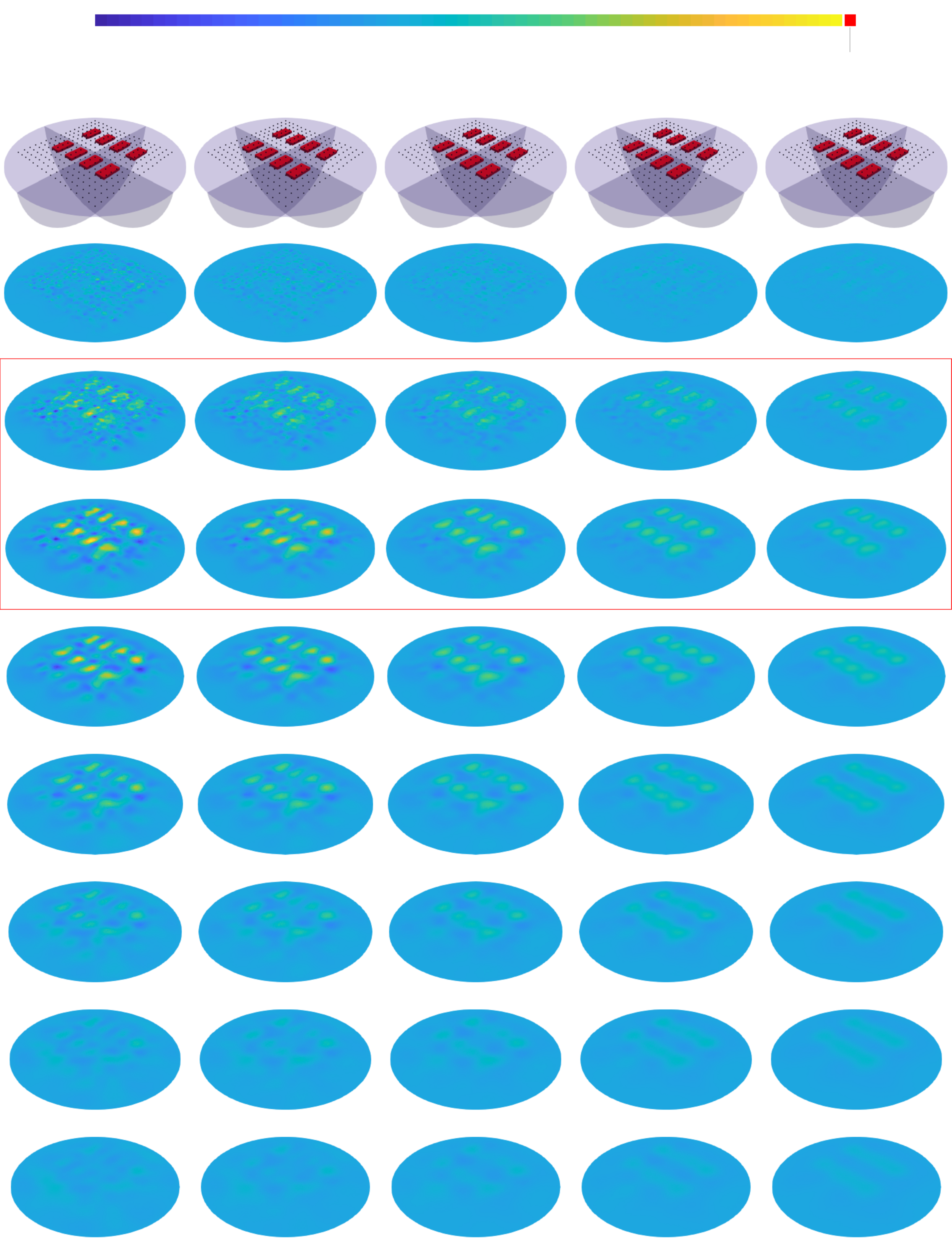}}%
    \put(0.1,1.28944155){\color[rgb]{0,0,0}\makebox(0,0)[t]{\lineheight{1.25}\smash{\begin{tabular}[t]{c}$2243$\end{tabular}}}}%
    \put(0.2960338,1.28944155){\color[rgb]{0,0,0}\makebox(0,0)[t]{\lineheight{1.25}\smash{\begin{tabular}[t]{c}$3046$\end{tabular}}}}%
    \put(0.49206761,1.28944155){\color[rgb]{0,0,0}\makebox(0,0)[t]{\lineheight{1.25}\smash{\begin{tabular}[t]{c}$3848$\end{tabular}}}}%
    \put(0.68810146,1.28944155){\color[rgb]{0,0,0}\makebox(0,0)[t]{\lineheight{1.25}\smash{\begin{tabular}[t]{c}$4651$\end{tabular}}}}%
    \put(0.88413527,1.28944155){\color[rgb]{0,0,0}\makebox(0,0)[t]{\lineheight{1.25}\smash{\begin{tabular}[t]{c}$5454$\end{tabular}}}}%
    \put(0.893238,1.23222438){\color[rgb]{0,0,0}\makebox(0,0)[t]{\lineheight{1.25}\smash{\begin{tabular}[t]{c}$7000$\end{tabular}}}}%
    \put(0.10013004,1.19139467){\color[rgb]{0,0,0}\makebox(0,0)[t]{\lineheight{1.25}\smash{\begin{tabular}[t]{c}$\beta = 10^{3.0}$\end{tabular}}}}%
    \put(0.29980495,1.19139467){\color[rgb]{0,0,0}\makebox(0,0)[t]{\lineheight{1.25}\smash{\begin{tabular}[t]{c}$\beta = 10^{3.5}$\end{tabular}}}}%
    \put(0.50013005,1.19139467){\color[rgb]{0,0,0}\makebox(0,0)[t]{\lineheight{1.25}\smash{\begin{tabular}[t]{c}$\beta = 10^{4.0}$\end{tabular}}}}%
    \put(0.69980497,1.19139467){\color[rgb]{0,0,0}\makebox(0,0)[t]{\lineheight{1.25}\smash{\begin{tabular}[t]{c}$\beta = 10^{4.5}$\end{tabular}}}}%
    \put(0.90013005,1.19139467){\color[rgb]{0,0,0}\makebox(0,0)[t]{\lineheight{1.25}\smash{\begin{tabular}[t]{c}$\beta = 10^{5.0}$\end{tabular}}}}%
    \put(0.00520156,0.93128739){\color[rgb]{0,0,0}\makebox(0,0)[lt]{\lineheight{1.25}\smash{\begin{tabular}[t]{l}$z = 0.00$\end{tabular}}}}%
    \put(0.20520156,0.93128739){\color[rgb]{0,0,0}\makebox(0,0)[lt]{\lineheight{1.25}\smash{\begin{tabular}[t]{l}$z = 0.00$\end{tabular}}}}%
    \put(0.40520157,0.93128739){\color[rgb]{0,0,0}\makebox(0,0)[lt]{\lineheight{1.25}\smash{\begin{tabular}[t]{l}$z = 0.00$\end{tabular}}}}%
    \put(0.60520155,0.93128739){\color[rgb]{0,0,0}\makebox(0,0)[lt]{\lineheight{1.25}\smash{\begin{tabular}[t]{l}$z = 0.00$\end{tabular}}}}%
    \put(0.80520158,0.93128739){\color[rgb]{0,0,0}\makebox(0,0)[lt]{\lineheight{1.25}\smash{\begin{tabular}[t]{l}$z = 0.00$\end{tabular}}}}%
    \put(0.00520156,0.79926787){\color[rgb]{0,0,0}\makebox(0,0)[lt]{\lineheight{1.25}\smash{\begin{tabular}[t]{l}$z = 2.50$\end{tabular}}}}%
    \put(0.20520156,0.79926787){\color[rgb]{0,0,0}\makebox(0,0)[lt]{\lineheight{1.25}\smash{\begin{tabular}[t]{l}$z = 2.50$\end{tabular}}}}%
    \put(0.40520157,0.79926787){\color[rgb]{0,0,0}\makebox(0,0)[lt]{\lineheight{1.25}\smash{\begin{tabular}[t]{l}$z = 2.50$\end{tabular}}}}%
    \put(0.60520155,0.79926787){\color[rgb]{0,0,0}\makebox(0,0)[lt]{\lineheight{1.25}\smash{\begin{tabular}[t]{l}$z = 2.50$\end{tabular}}}}%
    \put(0.80520158,0.79926787){\color[rgb]{0,0,0}\makebox(0,0)[lt]{\lineheight{1.25}\smash{\begin{tabular}[t]{l}$z = 2.50$\end{tabular}}}}%
    \put(0.00520156,0.66724903){\color[rgb]{0,0,0}\makebox(0,0)[lt]{\lineheight{1.25}\smash{\begin{tabular}[t]{l}$z = 5.00$\end{tabular}}}}%
    \put(0.20520156,0.66724903){\color[rgb]{0,0,0}\makebox(0,0)[lt]{\lineheight{1.25}\smash{\begin{tabular}[t]{l}$z = 5.00$\end{tabular}}}}%
    \put(0.40520157,0.66724903){\color[rgb]{0,0,0}\makebox(0,0)[lt]{\lineheight{1.25}\smash{\begin{tabular}[t]{l}$z = 5.00$\end{tabular}}}}%
    \put(0.60520155,0.66724903){\color[rgb]{0,0,0}\makebox(0,0)[lt]{\lineheight{1.25}\smash{\begin{tabular}[t]{l}$z = 5.00$\end{tabular}}}}%
    \put(0.80520158,0.66724903){\color[rgb]{0,0,0}\makebox(0,0)[lt]{\lineheight{1.25}\smash{\begin{tabular}[t]{l}$z = 5.00$\end{tabular}}}}%
    \put(0.00520156,0.5352276){\color[rgb]{0,0,0}\makebox(0,0)[lt]{\lineheight{1.25}\smash{\begin{tabular}[t]{l}$z = 7.50$\end{tabular}}}}%
    \put(0.20520156,0.5352276){\color[rgb]{0,0,0}\makebox(0,0)[lt]{\lineheight{1.25}\smash{\begin{tabular}[t]{l}$z = 7.50$\end{tabular}}}}%
    \put(0.40520157,0.5352276){\color[rgb]{0,0,0}\makebox(0,0)[lt]{\lineheight{1.25}\smash{\begin{tabular}[t]{l}$z = 7.50$\end{tabular}}}}%
    \put(0.60520155,0.5352276){\color[rgb]{0,0,0}\makebox(0,0)[lt]{\lineheight{1.25}\smash{\begin{tabular}[t]{l}$z = 7.50$\end{tabular}}}}%
    \put(0.80520158,0.5352276){\color[rgb]{0,0,0}\makebox(0,0)[lt]{\lineheight{1.25}\smash{\begin{tabular}[t]{l}$z = 7.50$\end{tabular}}}}%
    \put(0.00520156,0.40321197){\color[rgb]{0,0,0}\makebox(0,0)[lt]{\lineheight{1.25}\smash{\begin{tabular}[t]{l}$z = 10.00$\end{tabular}}}}%
    \put(0.20520156,0.40321197){\color[rgb]{0,0,0}\makebox(0,0)[lt]{\lineheight{1.25}\smash{\begin{tabular}[t]{l}$z = 10.00$\end{tabular}}}}%
    \put(0.40520157,0.40321197){\color[rgb]{0,0,0}\makebox(0,0)[lt]{\lineheight{1.25}\smash{\begin{tabular}[t]{l}$z = 10.00$\end{tabular}}}}%
    \put(0.60520155,0.40321197){\color[rgb]{0,0,0}\makebox(0,0)[lt]{\lineheight{1.25}\smash{\begin{tabular}[t]{l}$z = 10.00$\end{tabular}}}}%
    \put(0.80520158,0.40321197){\color[rgb]{0,0,0}\makebox(0,0)[lt]{\lineheight{1.25}\smash{\begin{tabular}[t]{l}$z = 10.00$\end{tabular}}}}%
    \put(0.00520156,0.27118983){\color[rgb]{0,0,0}\makebox(0,0)[lt]{\lineheight{1.25}\smash{\begin{tabular}[t]{l}$z = 12.50$\end{tabular}}}}%
    \put(0.20520156,0.27118983){\color[rgb]{0,0,0}\makebox(0,0)[lt]{\lineheight{1.25}\smash{\begin{tabular}[t]{l}$z = 12.50$\end{tabular}}}}%
    \put(0.40520157,0.27118983){\color[rgb]{0,0,0}\makebox(0,0)[lt]{\lineheight{1.25}\smash{\begin{tabular}[t]{l}$z = 12.50$\end{tabular}}}}%
    \put(0.60520155,0.27118983){\color[rgb]{0,0,0}\makebox(0,0)[lt]{\lineheight{1.25}\smash{\begin{tabular}[t]{l}$z = 12.50$\end{tabular}}}}%
    \put(0.80520158,0.27118983){\color[rgb]{0,0,0}\makebox(0,0)[lt]{\lineheight{1.25}\smash{\begin{tabular}[t]{l}$z = 12.50$\end{tabular}}}}%
    \put(0.00520156,0.13917428){\color[rgb]{0,0,0}\makebox(0,0)[lt]{\lineheight{1.25}\smash{\begin{tabular}[t]{l}$z = 15.00$\end{tabular}}}}%
    \put(0.20520156,0.13917428){\color[rgb]{0,0,0}\makebox(0,0)[lt]{\lineheight{1.25}\smash{\begin{tabular}[t]{l}$z = 15.00$\end{tabular}}}}%
    \put(0.40520157,0.13917428){\color[rgb]{0,0,0}\makebox(0,0)[lt]{\lineheight{1.25}\smash{\begin{tabular}[t]{l}$z = 15.00$\end{tabular}}}}%
    \put(0.60520155,0.13917428){\color[rgb]{0,0,0}\makebox(0,0)[lt]{\lineheight{1.25}\smash{\begin{tabular}[t]{l}$z = 15.00$\end{tabular}}}}%
    \put(0.80520158,0.13917428){\color[rgb]{0,0,0}\makebox(0,0)[lt]{\lineheight{1.25}\smash{\begin{tabular}[t]{l}$z = 15.00$\end{tabular}}}}%
    \put(0.00520156,0.00715215){\color[rgb]{0,0,0}\makebox(0,0)[lt]{\lineheight{1.25}\smash{\begin{tabular}[t]{l}$z = 17.50$\end{tabular}}}}%
    \put(0.20520156,0.00715215){\color[rgb]{0,0,0}\makebox(0,0)[lt]{\lineheight{1.25}\smash{\begin{tabular}[t]{l}$z = 17.50$\end{tabular}}}}%
    \put(0.40520157,0.00715215){\color[rgb]{0,0,0}\makebox(0,0)[lt]{\lineheight{1.25}\smash{\begin{tabular}[t]{l}$z = 17.50$\end{tabular}}}}%
    \put(0.60520155,0.00715215){\color[rgb]{0,0,0}\makebox(0,0)[lt]{\lineheight{1.25}\smash{\begin{tabular}[t]{l}$z = 17.50$\end{tabular}}}}%
    \put(0.80520158,0.00715215){\color[rgb]{0,0,0}\makebox(0,0)[lt]{\lineheight{1.25}\smash{\begin{tabular}[t]{l}$z = 17.50$\end{tabular}}}}%
  \end{picture}%
\endgroup%
   \unskip
  \caption{%
    A~fixed 3D example (third column from \cref{fig:ex3d},
    $N^\mathrm{ele}=289$, $N=452736$, $M=1564$)
    computed for series of regularization
    parameters $\beta$ using \cref{alg:gnstep-iterative}.
  }
  \label{fig:ex3d-beta}
\end{figure}
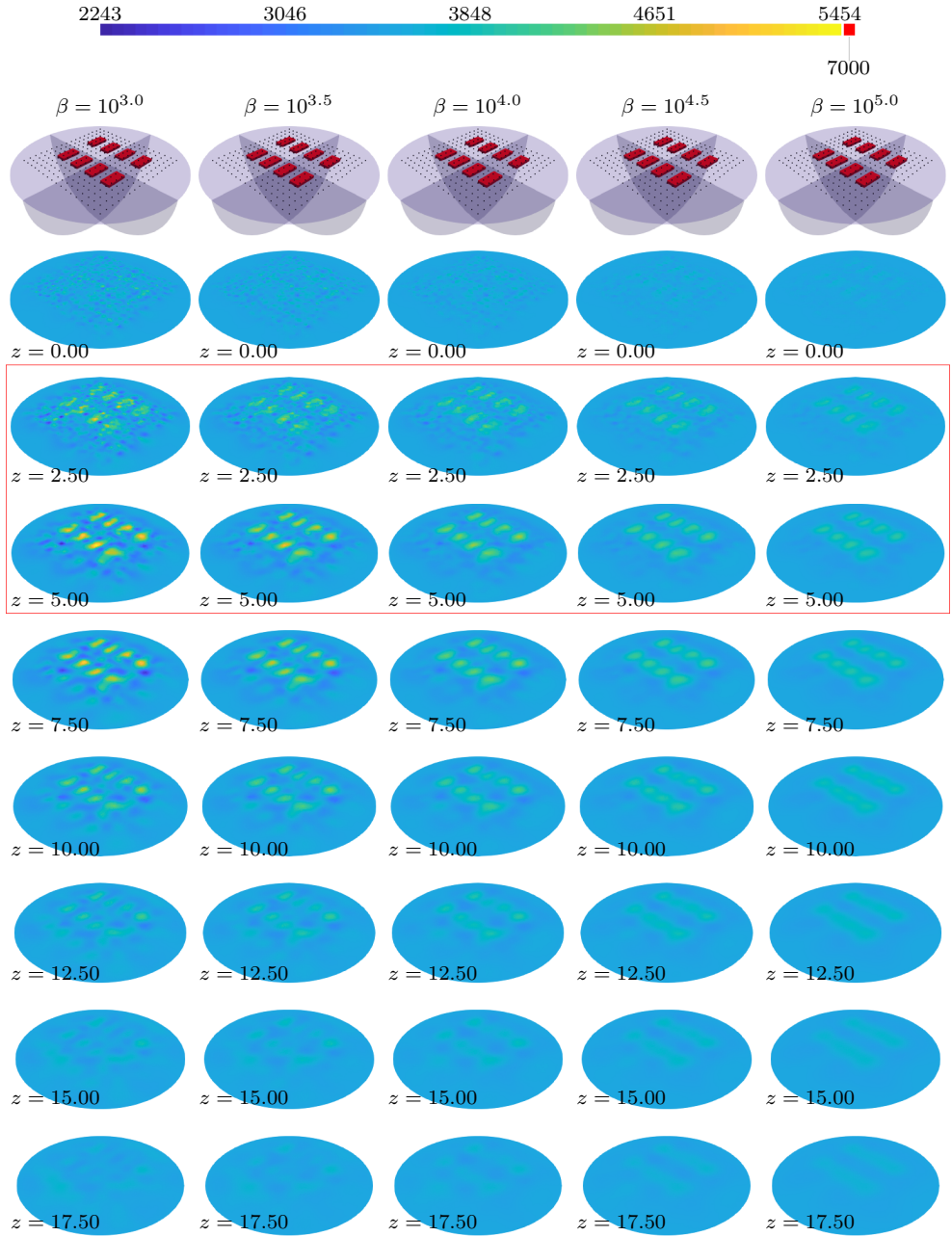

\begin{figure}[p]
  \setlength\figwidth{0.513\paperwidth}
  \setlength\figheight{0.132\paperwidth}
  \centering
  \footnotesize
  \capstart
  \mbox{%
%
\definecolor{mycolor1}{rgb}{0.00000,0.44700,0.74100}%
\definecolor{mycolor2}{rgb}{0.85000,0.32500,0.09800}%
\definecolor{mycolor3}{rgb}{0.92900,0.69400,0.12500}%
\definecolor{mycolor4}{rgb}{0.49400,0.18400,0.55600}%
\definecolor{mycolor5}{rgb}{0.46600,0.67400,0.18800}%
\definecolor{mycolor6}{rgb}{0.30100,0.74500,0.93300}%
\definecolor{mycolor7}{rgb}{0.63500,0.07800,0.18400}%
\begin{tikzpicture}

\begin{axis}[%
width=0.75\figwidth,
height=\figheight,
at={(0\figwidth,0\figheight)},
scale only axis,
xmode=log,
xmin=10000,
xmax=10000000000,
xminorticks=true,
xlabel style={font=\color{white!15!black}},
xlabel={$MN$},
ymode=log,
ymin=0.001,
ymax=10000,
yminorticks=true,
yticklabel style={rotate=90},
ylabel style={font=\color{white!15!black}},
ylabel={$t_{\mathrm{norm}}$ [s]},
axis background/.style={fill=white},
title style={font=\bfseries},
title={2D},
legend style={at={(1.03,1)}, anchor=north west, legend cell align=left, align=left, draw=white!15!black},
xlabel near ticks,
ylabel near ticks,
yticklabel style={anchor=south west,xshift=-2.5mm}
]
\addplot [color=mycolor1, only marks, mark=x, mark options={solid, mycolor1}]
  table[row sep=crcr]{%
38640	0.023274\\
224928	0.071311\\
1048760	0.190787\\
4316616	0.799432\\
17302984	3.242722\\
69155448	13.979481\\
273303712	63.390139\\
1096622704	330.029633\\
4371318032	1943.191575\\
};
\addlegendentry{$i=1$ (\cref{alg:gnstep-iterative})}

\addplot [color=mycolor2, only marks, mark=+, mark options={solid, mycolor2}]
  table[row sep=crcr]{%
38640	0.016647\\
224928	0.050751\\
1048760	0.224483\\
4316616	1.090607\\
17302984	3.883511\\
69155448	15.397441\\
273303712	69.60887\\
1096622704	357.531863\\
4371318032	1963.19775\\
};
\addlegendentry{$i=2$ (\cref{alg:gnstep-iterative})}

\addplot [color=mycolor3, only marks, mark=diamond, mark options={solid, mycolor3}]
  table[row sep=crcr]{%
38640	0.020769\\
224928	0.379772\\
1048760	6.46004\\
4316616	33.713602\\
17302984	180.168401\\
};
\addlegendentry{$i=1$ (\cref{alg:gnstep-iterative-nw})}

\addplot [color=mycolor4, only marks, mark=square, mark options={solid, mycolor4}]
  table[row sep=crcr]{%
38640	0.069421\\
224928	1.58141\\
1048760	39.690779\\
4316616	331.634947\\
17302984	2016.610792\\
};
\addlegendentry{$i=2$ (\cref{alg:gnstep-iterative-nw})}

\addplot [color=mycolor5, only marks, mark=triangle, mark options={solid, rotate=180, mycolor5}]
  table[row sep=crcr]{%
38640	0.043385\\
224928	0.252956\\
1048760	1.09116\\
4316616	4.663858\\
17302984	21.023429\\
69155448	107.233585\\
273303712	481.490465\\
};
\addlegendentry{$i=1$ (\cref{alg:gn-direct})}

\addplot [color=mycolor6, only marks, mark=triangle, mark options={solid, mycolor6}]
  table[row sep=crcr]{%
38640	0.040085\\
224928	0.228699\\
1048760	1.064458\\
4316616	4.605564\\
17302984	23.263214\\
69155448	105.507213\\
273303712	616.618842\\
};
\addlegendentry{$i=2$ (\cref{alg:gn-direct})}

\addplot [color=mycolor7, dashed]
  table[row sep=crcr]{%
38640	0.00347071343263912\\
4371318032	392.63955\\
};
\addlegendentry{slope $O(M N)$}

\addplot [color=mycolor1, dashdotted]
  table[row sep=crcr]{%
17302984	0.0941661317887352\\
69155448	0.765377646378474\\
273303712	6.0996216207825\\
1096622704	49.149895287164\\
4371318032	392.63955\\
};
\addlegendentry{slope $O(M^2 N)$}

\end{axis}
\end{tikzpicture}
}
  \rule[-0.5\baselineskip]{0pt}{\baselineskip}
  \mbox{%
%
\definecolor{mycolor1}{rgb}{0.00000,0.44700,0.74100}%
\definecolor{mycolor2}{rgb}{0.85000,0.32500,0.09800}%
\definecolor{mycolor3}{rgb}{0.92900,0.69400,0.12500}%
\definecolor{mycolor4}{rgb}{0.49400,0.18400,0.55600}%
\definecolor{mycolor5}{rgb}{0.46600,0.67400,0.18800}%
\definecolor{mycolor6}{rgb}{0.30100,0.74500,0.93300}%
\definecolor{mycolor7}{rgb}{0.63500,0.07800,0.18400}%
\begin{tikzpicture}

\begin{axis}[%
width=0.75\figwidth,
height=\figheight,
at={(0\figwidth,0\figheight)},
scale only axis,
xmode=log,
xmin=10000000,
xmax=10000000000,
xminorticks=true,
xlabel style={font=\color{white!15!black}},
xlabel={$MN$},
ymode=log,
ymin=10,
ymax=10000,
yminorticks=true,
yticklabel style={rotate=90},
ylabel style={font=\color{white!15!black}},
ylabel={$t_{\mathrm{norm}}$ [s]},
axis background/.style={fill=white},
title style={font=\bfseries},
title={3D \normalfont ($\beta=10^5$)},
legend style={at={(1.03,1)}, anchor=north west, legend cell align=left, align=left, draw=white!15!black},
xlabel near ticks,
ylabel near ticks,
yticklabel style={anchor=south west,xshift=-2.5mm}
]
\addplot [color=mycolor1, only marks, mark=x, mark options={solid, mycolor1}]
  table[row sep=crcr]{%
25961472	24.187135\\
191073792	139.583629\\
708079104	394.057708\\
1997130240	1208.23824\\
4405817600	2840.513461\\
};
\addlegendentry{$i=1$ (\cref{alg:gnstep-iterative})}

\addplot [color=mycolor2, only marks, mark=+, mark options={solid, mycolor2}]
  table[row sep=crcr]{%
25961472	20.979764\\
191073792	169.422769\\
708079104	392.27296\\
1997130240	1259.615544\\
4405817600	2988.638998\\
};
\addlegendentry{$i=2$ (\cref{alg:gnstep-iterative})}

\addplot [color=mycolor3, only marks, mark=diamond, mark options={solid, mycolor3}]
  table[row sep=crcr]{%
25961472	15.774641\\
191073792	94.939205\\
708079104	277.570395\\
1997130240	1363.21089\\
4405817600	2989.165872\\
};
\addlegendentry{$i=1$ (\cref{alg:gnstep-iterative-nw})}

\addplot [color=mycolor4, only marks, mark=square, mark options={solid, mycolor4}]
  table[row sep=crcr]{%
25961472	15.432866\\
191073792	72.149523\\
708079104	314.985516\\
1997130240	1203.991363\\
4405817600	3100.743216\\
};
\addlegendentry{$i=2$ (\cref{alg:gnstep-iterative-nw})}

\addplot [color=mycolor5, only marks, mark=triangle, mark options={solid, rotate=180, mycolor5}]
  table[row sep=crcr]{%
25961472	169.898276\\
191073792	1485.867739\\
708079104	7155.808383\\
};
\addlegendentry{$i=1$ (\cref{alg:gn-direct})}

\addplot [color=mycolor6, only marks, mark=triangle, mark options={solid, mycolor6}]
  table[row sep=crcr]{%
25961472	182.740012\\
191073792	1694.870818\\
708079104	6463.721538\\
};
\addlegendentry{$i=2$ (\cref{alg:gn-direct})}

\addplot [color=mycolor7, dashed]
  table[row sep=crcr]{%
25961472	12.3274797770293\\
4405817600	2092.0472986\\
};
\addlegendentry{slope $O(M N)$}

\addplot [color=mycolor1, dashdotted]
  table[row sep=crcr]{%
708079104	111.883399763784\\
1997130240	593.199595009187\\
4405817600	2092.0472986\\
};
\addlegendentry{slope $O(M^2 N)$}

\end{axis}
\end{tikzpicture}
}
  \vspace{-1\baselineskip}
  \caption{%
    Time for solving normal equations $t_\mathrm{norm}$
    for the 2D (top) and 3D (bottom) example.
    Each timing corresponds to runtime of
    \cref{alg:gnstep-iterative},
    \cref{alg:gnstep-iterative-nw},
    or \crefrange{ln:ldlt2}{ln:capa2} of \cref{alg:gn-direct}
    in each Gauss--Newton step~$i$.
  }
  \label{fig:timing}
  \vspace{1\baselineskip}
  \capstart
  \mbox{%
%
\definecolor{mycolor1}{rgb}{0.00000,0.44700,0.74100}%
\definecolor{mycolor2}{rgb}{0.85000,0.32500,0.09800}%
\definecolor{mycolor3}{rgb}{0.92900,0.69400,0.12500}%
\definecolor{mycolor4}{rgb}{0.49400,0.18400,0.55600}%
\begin{tikzpicture}

\begin{axis}[%
width=0.75\figwidth,
height=\figheight,
at={(0\figwidth,0\figheight)},
scale only axis,
xmode=log,
xmin=1000,
xmax=100000,
xminorticks=true,
xlabel style={font=\color{white!15!black}},
xlabel={$\beta$},
ymode=log,
ymin=243.242284,
ymax=2209.429871,
yminorticks=true,
yticklabel style={rotate=90},
ylabel style={font=\color{white!15!black}},
ylabel={$t_{\mathrm{norm}}$ [s]},
axis background/.style={fill=white},
title style={font=\bfseries},
title={3D \normalfont ($N^\mathrm{ele}=289$, $N=452736$, $M=1564$)},
legend style={at={(1.03,1)}, anchor=north west, legend cell align=left, align=left, draw=white!15!black},
xlabel near ticks,
ylabel near ticks,
yticklabel style={anchor=south west,xshift=-2.5mm}
]
\addplot [color=mycolor1, only marks, mark=x, mark options={solid, mycolor1}]
  table[row sep=crcr]{%
100000	365.867489\\
31622.7766016838	376.168628\\
10000	360.01807\\
3162.27766016838	349.972691\\
1000	342.117369\\
};
\addlegendentry{$i=1$ (\cref{alg:gnstep-iterative})}

\addplot [color=mycolor2, only marks, mark=+, mark options={solid, mycolor2}]
  table[row sep=crcr]{%
100000	381.490177\\
31622.7766016838	379.023704\\
10000	384.305232\\
3162.27766016838	487.99319\\
1000	361.357865\\
};
\addlegendentry{$i=2$ (\cref{alg:gnstep-iterative})}

\addplot [color=mycolor3, only marks, mark=diamond, mark options={solid, mycolor3}]
  table[row sep=crcr]{%
100000	243.242284\\
31622.7766016838	447.501069\\
10000	679.10105\\
3162.27766016838	1111.717576\\
1000	1883.311922\\
};
\addlegendentry{$i=1$ (\cref{alg:gnstep-iterative-nw})}

\addplot [color=mycolor4, only marks, mark=square, mark options={solid, mycolor4}]
  table[row sep=crcr]{%
100000	250.586372\\
31622.7766016838	400.377265\\
10000	738.284151\\
3162.27766016838	1243.808299\\
1000	2209.429871\\
};
\addlegendentry{$i=2$ (\cref{alg:gnstep-iterative-nw})}

\end{axis}
\end{tikzpicture}
}
  \rule[-0.5\baselineskip]{0pt}{\baselineskip}
  \mbox{%
%
\definecolor{mycolor1}{rgb}{0.00000,0.44700,0.74100}%
\definecolor{mycolor2}{rgb}{0.85000,0.32500,0.09800}%
\definecolor{mycolor3}{rgb}{0.92900,0.69400,0.12500}%
\definecolor{mycolor4}{rgb}{0.49400,0.18400,0.55600}%
\begin{tikzpicture}

\begin{axis}[%
width=0.75\figwidth,
height=\figheight,
at={(0\figwidth,0\figheight)},
scale only axis,
xmode=log,
xmin=1000,
xmax=100000,
xminorticks=true,
xlabel style={font=\color{white!15!black}},
xlabel={$\beta$},
ymode=log,
ymin=100,
ymax=2739,
yminorticks=true,
yticklabel style={rotate=90},
ylabel style={font=\color{white!15!black}},
ylabel={$n_{\mathrm{iter}}$},
axis background/.style={fill=white},
title style={font=\bfseries},
title={3D \normalfont ($N^\mathrm{ele}=289$, $N=452736$, $M=1564$)},
legend style={at={(1.03,1)}, anchor=north west, legend cell align=left, align=left, draw=white!15!black},
xlabel near ticks,
ylabel near ticks,
yticklabel style={anchor=south west,xshift=-2.5mm}
]
\addplot [color=mycolor1, only marks, mark=x, mark options={solid, mycolor1}]
  table[row sep=crcr]{%
100000	124\\
31622.7766016838	121\\
10000	116\\
3162.27766016838	109\\
1000	103\\
};
\addlegendentry{$i=1$ (\cref{alg:gnstep-iterative})}

\addplot [color=mycolor2, only marks, mark=+, mark options={solid, mycolor2}]
  table[row sep=crcr]{%
100000	126\\
31622.7766016838	129\\
10000	130\\
3162.27766016838	130\\
1000	124\\
};
\addlegendentry{$i=2$ (\cref{alg:gnstep-iterative})}

\addplot [color=mycolor3, only marks, mark=diamond, mark options={solid, mycolor3}]
  table[row sep=crcr]{%
100000	291\\
31622.7766016838	487\\
10000	818\\
3162.27766016838	1363\\
1000	2228\\
};
\addlegendentry{$i=1$ (\cref{alg:gnstep-iterative-nw})}

\addplot [color=mycolor4, only marks, mark=square, mark options={solid, mycolor4}]
  table[row sep=crcr]{%
100000	297\\
31622.7766016838	520\\
10000	921\\
3162.27766016838	1629\\
1000	2739\\
};
\addlegendentry{$i=2$ (\cref{alg:gnstep-iterative-nw})}

\end{axis}
\end{tikzpicture}
}
  \vspace{-1\baselineskip}
  \caption{%
    Time for solving normal equations $t_\mathrm{norm}$ (top)
    and number of MINRES iterations (bottom)
    for a~fixed 3D example (third column from \cref{fig:ex3d})
    with series of regularization parameters~$\beta$.
  }
  \label{fig:timing-beta}
\end{figure}

\medskip\paragraph{2D test case}
We consider the half-disk domain $\Omega\coloneqq\{(x,z)\in{\mathbb R}^2,\,z>0,\,\allowbreak\sqrt{x^2+z^2}<80\}$.
The line $\{z=0\}$ represents the ground surface where measurements are taken using electrodes placed as described
in \cref{fig:survey}.
Following geophysical convention, the half-space $\{z>0\}$ represents the subsurface consisting of a~medium with (here a~priori known) conductivity distribution~$\sigma_\mathrm{true}$ as in \cref{fig:ex2d} on the left displaying
a~series of increasingly finer anomalous conductivity patterns imposed on a~background medium of constant conductivity.
The opposite side $\{z<0\}$ represents the air half-space of negligible conductivity, which is thus excluded from the domain and modeled by a~vanishing normal component of the electric field~\cref{eq_model_gn} on ${\gamma_{\mathrm N}}\coloneqq\{z=0\}$.
For simplicity we consider~\cref{eq_model_gd} on ${\gamma_{\mathrm D}}\coloneqq{\partial\Omega}\setminus\overline{\gamma_{\mathrm N}}$.
This description fully specifies the functions
$m\mapsto g^i(m)$ and $m\mapsto J^i(m)$, $i=1,2,\ldots,M$.

For the configurations, the first five of which are indicated in \cref{fig:ex2d}, we compute the finite element approximations of the quantities ${g^i_{\mathrm{obs}}}\coloneqq g^i(\log\sigma_\mathrm{true})$, $i=1,2,\ldots,M$, which serve as the (synthetic) observational data for inversion.
Note that this data is noisy due to the discretization error (although the meshes used to generate the values ${g^i_{\mathrm{obs}}}$ are finer compared to the meshes for the inversion).
The reference value is taken to be ${m_\mathrm{ref}}\coloneqq\log\sigma_\mathrm{ref}$ as in
\cref{fig:ex2d} and ${\Gamma_{\mathrm D}}$ in~\cref{eq:ls} is taken as
${\Gamma_{\mathrm D}}\coloneqq{\partial\Omega}$. Two Gauss--Newton steps with a~fixed value
of the regularization parameter $\beta\coloneqq 0.1$
are performed and the resistivity distributions
in \cref{fig:ex2d} (on the right) are obtained.
The meshes for the inversion (see \cref{fig:ex2d})
are a~priori refined around the electrode positions,
which are at the surface $\{z=0\}$, so that the meshes
scale as $N=O(N^\mathrm{ele})$; see \cref{tab:experiments}.

\medskip\paragraph{3D test case}
Here we consider the semi-spherical domain $\Omega\coloneqq\{(x,y,z)\in{\mathbb R}^3,\,z>0,\,\allowbreak\sqrt{x^2+y^2+z^2}<80\}$.
The measurements are again taken on the surface $\{z=0\}$ using the grid of electrodes shown in \cref{fig:ex3d}.
One uses the pole-dipole scheme (as described in \cref{fig:survey}) along the $x$-direction for all possible $y=\mathrm{const}$ profiles and then the same in the $y$-direction
for all possible $x=\mathrm{const}$ profiles.
By analogy, the true resistivity model is also constructed in a~Cartesian product fashion;
see \cref{fig:ex3d}.

The remaining details are analogous to the 2D test case above with the exception that different values of the regularization parameter $\beta$ were necessary to obtain good reconstructions.
The question of choosing the best value of the regularization parameter are beyond the scope of this paper.
Nevertheless we experimented with a~number of choices and noticed how this affects the performance of the  algorithms.
For the sake of illustration, we indicate in \cref{fig:ex3d-beta} the effect of the regularization parameter on the reconstructed conductivity.

\afterpage{\clearpage}

\medskip
The linear systems resulting from 2D discretizations of~\eqref{eq_model} were solved using a~sparse direct method and in 3D using conjugate gradient iteration preconditioned by an algebraic multigrid cycle.
These choices make the approximation of $A^{-1}_{\exp(m)}$ in~\eqref{eq:ert_qoi2}, and in turn computation of $\bm{g}_{\bm{m}}$ and ${\boldsymbol{J}}_{\bm{m}}$, sufficiently inexpensive and scalable, leaving the main effort in the solution of~\eqref{eq:block_mixed_h}, which is the primary concern of this work.

The numerical experiments were implemented using Matlab,
HSL\_MI20~\cite{hslmi20}, and Gmsh~\cite{gmsh}. The plots
were produced using matlab2tikz~\cite{matlab2tikz}
and PyVista~\cite{sullivan2019pyvista}.
The complete code for reproducing the experiments
is available as~\cite{paper-woodbury-code}.

\Cref{tab:experiments} and \cref{fig:timing} show that the computational cost of the examples agrees with the expected complexity as predicted in \cref{sec:iterative}.
In particular, we can see that the dominating cost of \cref{alg:gnstep-iterative} is $O(MN)$ but we can see the
$O(M^2 N)$ term becoming effective for larger values of~$M$.
The Cholesky factorization for $O(M^3)$ (value $t_{\labelcref{ln:chol}}$ in \cref{tab:experiments})
and the matrix-matrix product for $O(M^2 N)$ (value $t_{\labelcref{ln:C}}$ in \cref{tab:experiments}) have a~small multiplicative constant as these would typically run very efficiently in LAPACK and BLAS, respectively.
Nevertheless, it is clear that $O(M^2 N)$ will dominate for larger problems.

\Cref{fig:timing-beta} shows the performance of a~fixed 3D test case depending on the value of the regularization parameter~$\beta$.
In particular, \cref{alg:gnstep-iterative} is seen to exhibit robust performance independent of~$\beta$.
For this one has to pay the price of computing and factoring the capacitance matrix.
\cref{alg:gnstep-iterative-nw}, on the other hand, shows strong dependence of the required number of MINRES steps
on the value of~$\beta$ and $MN$.
Although \cref{alg:gnstep-iterative-nw} may sometimes be a~less expensive alternative, \cref{alg:gnstep-iterative}
should generally be preferred for its robustness.
To this end we also note, that with a~better implementation of the solver for $\boldsymbol{\hat S}^{-1}$, one might achieve more favorable timings for the computation of the capacitance matrix.
We have used HSL\_MI20~\cite{hslmi20}, which is fully sequential in contrast to the threaded BLAS used in other parts of the code; 8~threads were used where applicable.
Moreover, HSL\_MI20 only implements $\boldsymbol{\hat S}^{-1}\bm{z}$ for a~single-column vector~$\bm{z}$, but we need, on \cref{ln:H} in \cref{alg:gnstep-iterative}, to apply $\boldsymbol{\hat S}^{-1}$ to all the $M$~columns of~${\boldsymbol{J}}_{\bm{m}}^\top$.
This operation therefore runs sequentially column-by-column, which is certainly not optimal in utilizing theoretical floating-point performance and memory bandwidth of the machine.
This implementation drawback penalizes \cref{alg:gnstep-iterative} in this experimental performance assessment and it should be kept on mind that \cref{alg:gnstep-iterative} can be more favorable than \cref{alg:gnstep-iterative-nw} whenever a~suitable AMG implementation is available.

Furthermore, we have observed that, for lower values of the regularization parameter~$\beta$ (for example, the 2D case with $\beta=0.001$, which is not shown in the paper), the solutions produced with preconditioners~$\boldsymbol{\hat P}_{\beta,\bm{m}}^{-1}$ and~$\boldsymbol{\hat P}^{-1}$ may differ significantly although they were solved to the same residual accuracy in the Euclidean norm:
$\norm{\bm{r}_k}_2/\norm{\bm{b}}_2\leq10^{-7}$.
Note that this has always been used as the stopping criterion in MINRES although the minimization intrinsic to the preconditioned MINRES process minimizes a~different quantity;\footnote{%
This is Matlab's actual behavior:
  \texttt{MINRES(A, b, tol, maxit, M1, M2, x0)}
  mathematically means, for $\boldsymbol{A}$ symmetric and
  $\boldsymbol{P}$ symmetric positive definite,
  $
    \norm{\bm{r}_k}_{\boldsymbol{P}^{-1}}
    = \min_{p\in\mathcal{P}^0_k}
    \norm{p(\boldsymbol{A}\boldsymbol{P}^{-1})\bm{r}_0}_{\boldsymbol{P}^{-1}}
  $,
  $
    \norm{\bm{R}_k}_{\boldsymbol{P}}
    = \min_{p\in\mathcal{P}^0_k}
    \norm{p(\boldsymbol{P}^{-1}\boldsymbol{A})\bm{R}_0}_{\boldsymbol{P}}
  $, or
  $
    \norm{\boldsymbol{\rho}_k}_2
    = \min_{p\in\mathcal{P}^0_k}
    \norm{p(\boldsymbol{L}^{-1}\boldsymbol{A}\boldsymbol{L}^{-\top})\boldsymbol{\rho}_0}_2
  $,
  where all of these are equivalent formulations through
  $\bm{r}_k=\bm{b}-\boldsymbol{A}\bm{x}_k$,
  $\bm{R}_k=\boldsymbol{P}^{-1}\bm{r}_k$,
  $\boldsymbol{\rho}_k=\boldsymbol{L}^{-1}\bm{r}_k$,
  $\boldsymbol{L}\boldsymbol{L}^\top=\boldsymbol{P}$,
  and preconditioner~$\boldsymbol{P}$ is given by \texttt{M1} and
  \texttt{M2} as per the function's docstring.
  On the other hand, the function uses the Euclidean
  stopping criterion
  $\norm{\bm{r}_k}_2 \leq \texttt{tol}\, \norm{\bm{b}}_2$,
  regardless of the preconditioner and the initial guess.
} see~\cref{eq:minres}.

\section{Conclusion and outlook} \label{sec:outro}

We have formulated a~nonlinear parameter identification problem subject to $H^1$ regularization and its Gauss--Newton linearization as a~second-order boundary value problem including a~consistent interpretation of possible choices of boundary conditions as they result from the nature of the regularization procedure.
For a~standard inf-sup stable mixed discretization, we have proposed a~number of efficient and robust solution strategies of the linear systems arising from the Gauss--Newton linearization.
The proposed methods included a~direct method, a~preconditioned iterative scheme based on the Woodbury formula, and a~preconditioned iterative scheme in which the low-rank perturbation is not accounted for by the preconditioner and must be compensated by the Krylov iteration.
In a~series of extensive numerical experiments, we have performed scaling tests w.r.t.\ the relevant problem parameters for a~challenging parameter identification problem arising in electrical resistivity tomography.

In \cref{tab:algocomp} we summarize our findings concerning the interplay between efficiency and quality of the two considered preconditioners.
Note that, regarding the indicated scaling of MINRES iterations required to solve the linear system to prescribed tolerance, we do not have a~rigorous theoretical argument, but merely empirical findings specific to the class of problems we solved; see \cref{tab:experiments} for the observed dependence on~$M$ and \cref{fig:timing-beta} for the dependence on~$\beta$.
Specifically, the MINRES convergence behavior for~$\boldsymbol{\hat P}^{-1}$ is sure to be problem dependent and likely depends on the distribution of singular values of~${\boldsymbol{J}}_{\bm{m}}$.
The simple parametrization $M^{\gamma_1}$ observed here may only apply when ${\boldsymbol{J}}_{\bm{m}}$'s are selected from a~narrow class.
We have observed that the variant with the full preconditioner~\smash{$\boldsymbol{\hat P}_{\beta,\bm{m}}^{-1}$} exhibits robustness of convergence for a~range of parameter values $M$, $N$, and~$\beta$,
while the cheaper preconditioner~$\boldsymbol{\hat P}^{-1}$ can suffer from slow convergence (\cref{fig:timing-beta}) or
even stagnation (\cref{tab:experiments}).
\begin{table}
  \capstart
  \setlength\belowcaptionskip{0pt}  
  \caption{}
  \footnotesize
  \centering
  \begin{tabular}{p{0.35\linewidth}|cc}
    &
    \cref{alg:gnstep-iterative}
    &
    \cref{alg:gnstep-iterative-nw}
    \\
    \hline
    employed preconditioner
    &
    \rule{0pt}{2.5ex}%
    $\boldsymbol{\hat P}_{\beta,\bm{m}}^{-1}$
    &
    \rule{0pt}{2.5ex}%
    $\boldsymbol{\hat P}^{-1}$
    \\
    ${\boldsymbol{J}}_{\bm{m}}^\top{\boldsymbol{J}}_{\bm{m}}$ handled by
    &
    Woodbury
    &
    Krylov
    \\
    cost per MINRES iteration
    &
    $O(M^2N)$
    &
    $O(MN)$
    \\
    number of MINRES iterations
    &
    $O(1)$
    &
    $O(M^{\gamma_1} \beta^{-\gamma_2})$, $\gamma_1$, $\gamma_2>0$
    \\
    overall robustness
    &
    \newcheckmark
    &
    \newcrossmark
  \end{tabular}
  \label{tab:algocomp}
\end{table}

In future work we would like to investigate data sparse approximation and fast solution of the capacitance matrix equation in applying the Woodbury formula, e.g., using $\mathcal{H}$-matrix methods and/or randomized low-rank approximations.
This would allow applying the preconditioner $\boldsymbol{\hat P}_{\beta,\bm{m}}^{-1}$ with a~lower complexity than $O(M^2N)$.

\section*{Acknowledgment}
The authors are grateful to Mathias Scheunert (Technische Universit\"at Bergakademie Freiberg) for programming mesh generation for the computational examples.

\begingroup
\raggedbottom
\interlinepenalty=10000

\bibliographystyle{siamplain}
\endgroup

\end{document}